\documentclass[12pt]{article}
\usepackage[english]{babel}
\usepackage{a4wide}
\usepackage{amssymb}
\usepackage{amsmath}
\usepackage{mathdots}
\usepackage{theorem}
\usepackage{exon}
\usepackage[latin1]{inputenc}
\usepackage[T1]{fontenc}

\title{A numerical proof of the Gr{\"u}nbaum conjecture}
\author{David Hermann \footnote{david.hermann@imj-prg.fr}} 

\begin{document}

\maketitle

\begin{abstract}
  Chalmers and Lewicki published in 2010 a very technical proof
  of the Gr{\"u}nbaum conjecture. Here is a simpler proof,
  partially based on numerical arguments.
\end{abstract}

The Hahn-Banach theorem states that onto each line in every normed space,
there is a unitary projection, and Kadec and Snobar \cite{ks}
proved (using John's ellipsoid) 
that onto each $n$-dimensional subspace of any real normed space,
there is a projection with norm at most $\lambda_n \leq \sqrt{n}$\,.
Gr{\"u}nbaum \cite{g} conjectured  that 
$\lambda_2=4/3=1.333...<1.414...=\sqrt{2}$\,, 
which is the {\sl projection constant} of the plane of equation $x_1+x_2+x_3=0$ in 
$(\, \MR^3 \vge \norme{\,\,}_\infty )$
whose norm is hexagonal, hence the $4/3$...
Several attempts have been made to prove this conjecture:
K{\"o}nig and Tomczak-Jaegermann published in \cite{ktj} a proof that was shown
incomplete by Chalmers and Lewicki, who gave their own (a bit intricate) proof in \cite{cl}. 
Here is a simpler proof, mostly based on their works, 
and partially on a few numerical studies of extrema of functions of 3 variables. 
Using arguments due to Lewis \cite{l},
K{\"o}nig and Tomczak-Jaegermann proved that
if $O^{N}_{\,n}$ denotes
the space of $n$ orthonormal vectors in the standard 
Euclidean space $\MR^N$ for every integers $1\leq n \leq N$ and if we set 
for all $(\,u\vge (x \vge y))\in O^{N}_{\,1} \times O^{N}_{\,2}$\,:
\[
\Phi_2^{N}\, (u\vge  (x\vge y)) \,\,\, := \,\,\,
\sum_{1\,\leq\,\,i\vge j\,\,\leq \,N} u_i\,\, u_j\,\,
\Absole{  x_i \, x_j\,+\, y_i \, y_j  }
\quad ,
\]
it suffices to prove that for each integer $N\geq 3$ we have:
\[
\lambda_{\,2}^N 
\,\,\,:=\,\,\,
\max_{O^{N}_{\,1} \times \, O^{N}_{\,2}} \Phi_2^{N}
 \,\,\,\leq  \,\,\, \frac{\,4\,}{3} \quad ,
\]
and here we leave the realm of Banach spaces geometry,
and from now on our only goal will be to estimate the maximum
of this function.
\medskip

For any $(x\vge y)\in O^N_2$\,, the symmetric matrix 
$P_{x\vge y} := ((x_i \, x_j\,+\, y_i \, y_j))_{i\vge j}\in Sym_N$ represents the orthogonal projection onto the plane $\calp_{x\vge y}\in O_2^N/O_2$ generated by $x$ and $y$ in the euclidean space $\MR^N$,  that is then regarded  as sitting into 
$(\, \MR^N \vge \norme{\,\,}_\infty )$: this, and a clever use of the Cauchy-Schwartz inequality in a probabilistic approach explain the formulae above.
Moreover, it explains the invariance of $\Phi_2^N$ under the right action of $O_2$ and the left actions of $O_N$.
Our proof runs as follows: the sequence $(\lambda_2^N )_N$ clearly increases, and if 
we have
$\lambda_{\,2}^{N}=\Phi_{N}\, (u\vge (x\vge y)) $ where
$x_N=y_N=0$\,, we
easily get $\lambda_{\,2}^{N} \leq \lambda_{\,2}^{N-1}$ and we
can conclude by induction (we will compute $\lambda_{\,2}^{\,3}=4/3$ 
in section \ref{n=3}), and else in
sections \ref{symetries} to \ref{inversible}
we will reduce by symmetry to the case where $N=2\,s+1$
for an integer $s\geq 2$  and:
\[
\Phi_2^{N}\, (u\vge (x\vge y)) 
\,\,\, = \,\,\,
\phi_A \,(u\vge (x\vge y)) 
\,\,\, = \,\,\,
\sum_{1\,\leq\,\,i\vge j\,\,\leq \,N} a_{\, i\vge j}\,\,  
u_{\, i}\,\,  u_j\,\,
\big(\, x_i \, x_j\,+\, y_i \, y_j \,\big) \quad ,
\]
where $A=((a_{\,i\vge j}))_{\,1\leq\, i\vge j\, \leq\, N}$
is a given symmetric matrix with coefficients in 
$\{-1\vge 1\}$\,.
In section \ref{essentiels}, the study of the critical
points of $\phi_A$ will allow us to prove that if
$\lambda_{\,2}^{N}>\lambda_{\,2}^{N-1}$, then there exists two real numbers
$\alpha$ and $\beta$ such that $1/3 < \beta \leq \alpha \leq 1$ and 
$\alpha+\beta > 4/3$\,, and a finite sequence:
\, $
0 \leq
\theta_1 < \theta_2 < \cdots  < \theta_s  < \theta_{s+1}  
< \theta_1+\pi < 2\,\pi
$ \,
satisfying the induction relation for each $1\leq k \leq s$\,:
\[
\sin \Big(\,\frac{\theta_{k+1}-\theta_{k}}{2}\,\Big)
=\frac{1}{2\,s+1}\,\sqrt{\,
\frac{1}{\beta^2}
-\Big(\,\frac{1}{\beta^2}-\frac{1}{\alpha^2}\,\Big)\,
\sin^2\Big(\,\frac{\theta_{k}+\theta_{k+1}}{2}\,\Big)\,}
\quad ,
\]
the boundary condition:
\, $
\dsm 
\beta\,   \sin \theta_{1}\, \sin \theta_{s+1}  +
\alpha\,   \cos \theta_1  \,  \cos \theta_{s+1} 
=\frac{1}{2s+1}
$ \,
and the equation:
\[
\frac{1}{2\,s+1}\,\,\sum_{k=1}^{s} \,\, 
\cos\big(\,\theta_k+\theta_{k+1}\,\big) 
\,=\,
\frac{(s+1)\,\beta-s\,\alpha}{(2s+1)\,(\alpha+\beta)} 
-\beta\,\sin \theta_{1}\, \sin \theta_{s+1}  
\quad .
\]
Then we will prove numerically that such a sequence doesn't exist,
which proves Gr{\"u}nbaum conjecture by induction: for  $s$ large enough,
each $\theta_k$ approximates 
$\dsm y\Big(\frac{2\,k-2}{2\,s+1}\Big)$\,, 
where the function $y\dpe [\,0\vge 1\,] \dans \MR$ satisfies the
differential equation: 
\[
y'=\frac{1}{\beta}\,\sqrt{\,1
-\Big(1-\frac{\beta^2}{\alpha^2}\,\Big)\,\sin^2 y\,}
\]
with the initial condition
$y(0)=\theta_1\in [\,0 \vge \pi\,[$\,,  
the boundary condition: 
\[
y(0) \leq y(1) \leq y(0)+ \pi
\quad \hbox{and} \quad
\alpha\,   \cos y(0)  \,  \cos y(1)
+ \beta\,   \sin y(0)\, \sin  y(1)
\,=\,
0
\]
and the integral equation 
(which is indeed another boundary condition):
\[
\int_{y(0)}^{\,y(1)}
\frac{\alpha\,\beta \,\cos (2\,x)\,\,dx}
{\sqrt{\alpha^2\, \cos^2 x+\beta^2\, \sin^2  x}}
+\frac{\alpha-\beta}{\alpha+\beta} 
+\frac{\alpha\,\beta\,\sin (\, 2\,y(0)\,) }
{\sqrt{\rule{0cm}{0.36cm}\alpha^2\, \cos^2 y(0)+\beta^2\, \sin^2  y(0)}}
\,=\,
0
\]
and we will get in sections \ref{equadiff} to \ref{en}
explicit estimates for this kind of ``middle-point at the goal method'', 
which will allow us to reduce the problem to the numerical study of
the minimum of a function of $3$ variables, which will be done
(using Maple) in section \ref{numerique-ed},
where we will conclude for $s\geq 15$ since these conditions are
incompatible. 
The remaining cases where $2\leq s \leq 14$ will be treated similarly
in section \ref{iteration}, 
but will require  procedures in C (using Code::Blocks)
in order to keep the computation time reasonable: on my own PC,
the Maple procedures take less than $2$ hours and the C procedures take
less than $6$ hours.
Finally, in the last sections we will get the estimates on the partial derivatives
of the relevant functions that are needed in sections  \ref{numerique-ed}
and \ref{iteration} in order to conclude.

\section{The symmetries of $\Phi_N$ and the matrix of signs}
\label{symetries}

For each integer $N\geq 2\,,$ set as above:
\quad $\dsm \lambda_{\,2}^N \, = \,
\max_{O^{N}_{\,1} \times O^{N}_{\,2}}
\Phi_{\,2}^{N} $ \quad where:
\[
\Phi_{\,2}^{N}\, (u\vge (x\vge y)) \, = \,
\sum_{1\,\leq\,\,i\vge j\,\,\leq \,N} u_i\,\, u_j\,\,
\Absole{  x_i \, x_j\,+\, y_i \, y_j  } \quad \hbox{ for all }
(u\vge (x \vge y))\in O^{N}_{\,1} \times O^{N}_{\,2} \quad ,
\]
thus the sequence $(\lambda_{\,2}^N)_{\,N\geq 2}$ is increasing
and \,
$\Phi_{\,2}^{N}\, (u\vge (x\vge y)) \,\leq \,
\Phi_{\,2}^{N}\, (\, \absol{u}\vge (x\vge y))$
\, where:
\[
\absol{(u_1\vge \dots \vge u_N)}\,=\,
(\,\absol{u_1}\vge \dots \vge \absol{u_N}\,) \in O^{N}_{\,1} \quad ,
\]
thus the maximum of $\Phi_{\,2}^N$ is attained at a point where each $u_i$ is
nonnegative. We can permute all the coordinates without changing 
$O^{N}_{\,1} \times O^{N}_{\,2}$ and the value of $\Phi_2^N$\,, and for all 
$\eps \in \{-1\vge 1\}^N$
let $s_\eps \dpe (v_1\vge \dots \vge v_N)\mapsto
(\eps_1\, v_1\vge \dots \vge \eps_N\,v_N)$\,, so
$s_\eps \times s_\eps$ preserves $O^{N}_{\,2}$ and we have:
$\Phi_{\,2}^N (u \vge (s_\eps(x)\vge s_\eps(y)))=
\Phi_{\,2}^{N}\, (u\vge (x\vge y))$ for all
$(u\vge (x \vge y))\in O^{N}_{\,1} \times O^{N}_{\,2}$\,, 
thus we are free to change the sign of $(x_i\vge y_i)$ for each
$1\leq i \leq N$\,.
If $x_N=y_N=0$\,, let $x=(\,x^\ast \vge 0)$\,, $y=(\,y^\ast \vge 0)$
and $u=(\, u^\ast \vge u_N)$ where
$u^\ast \vge x^\ast \vge y^\ast \in \MR^{N-1}$ 
and $u_N\in [-1\vge 1]$\,, thus
$(\,x^\ast \vge y^\ast)\in O^{N}_{\,2}$ and 
$u^\ast=\sqrt{1-u_N^2\,}\,\, u^\dagger$ where 
$u^\dagger \in  O^{N}_{\,1}$\,, and we get: 
\[
\Phi_{\,2}^{N}\, (u\vge (x\vge y)) \, = \,
\Phi_{N-1}\, (u^\ast\vge (x^\ast\vge y^\ast)) \, = \,
\big(\,1-u_N^2\big)\,\Phi_{N-1}\, (u^\dagger\vge (x^\ast\vge y^\ast)) 
\, \leq \, \lambda_{\,2}^{N-1} \quad ,
\]
thus, arguing by induction, we can assume that
$\Phi_{\,2}^N$ attains its maximum at a point where
$(\,x_i\vge y_i)\neq 0$ for all $1\leq i \leq N$\,.
For all $(\,x \vge y)\in  O^{N}_{\,2}$\,, 
the matrix
$((x_i\,x_j+y_i\,y_j))_{\,1\leq\, i\vge j\, \leq N }$
represents the orthogonal projection onto the plane 
$\calp_{x\vge y}$ generated by $(x \vge y)$\,, thus if  
$(\,\widetilde{x} \vge \widetilde{y}\,)$
is another orthonormal basis of $\calp_{x\vge y}$ we have:
$\Phi_2^{N}\, (\,u\vge (\widetilde{x}\vge \widetilde{y}\,)) \, = \,
\Phi_2^ {N}\, (u\vge (x\vge y))$\,.
If $N=2$\,, we can thus assume that $(\,x \vge y)$ is
the canonical basis to get $\lambda_{\,2}^{\,2}=1$\,,
which is geometrically obvious. 
\bigskip

For each $(\,x \vge y)\in  O^{N}_{\,2}$ and every
$1\,\leq\,\,i\vge j\,\,\leq \,N $\,, let 
$a_{\,i\vge j}\, (x\vge y)\, \in \,  \{-1\vge +1\}$ 
be the sign of $x_i \, x_j\,+\, y_i \, y_j$
where $0$ is positive, thus 
$A\,(x\vge y)\in \cala_N \subset Sym (N) \subset M_N (\MR)$ 
is a symmetric matrix 
having coefficients in $\{-1\vge 1\}$ and diagonal coefficients equal
to $1$\,. For all $A\in \cala_N$ and 
$(u\vge (x \vge y))\in O^{N}_{\,1} \times O^{N}_{\,2}$\,, we have:
\[
\phi_A \,(u\vge (x\vge y)) 
\,\,\, = \,\,\,
\sum_{1\,\leq\,\,i\vge j\,\,\leq \,N} a_{\, i\vge j}\,\,  
u_{\, i}\,\,  u_j\,\,
\big(\, x_i \, x_j\,+\, y_i \, y_j \, \big)
\,\,\, \leq \,\,\,
\Phi_{\,2}^N \,(\,\absol{u}\vge (x\vge y)) 
\]
and 
$\Phi_{\,2}^N(u\vge (x \vge y))=\phi_{A(x\vge y)} (u \vge (x \vge y)) $\,,
thus if $\Phi_{\,2}^N$ attains its maximum at $(u\vge (x \vge y))$\,,
it is also the maximum of the differentiable function $\phi_A$ where 
$A=A(x\vge y)$ depends only on $\calp_{x\vge y}$\,: 
this trick is due to Chalmers and Lewicki, 
like the beginning of the next section.

\section{The critical points of $\phi_A$}\label{extrema}

For each integer $N\geq 2$ and every symmetric matrix $B\in Sym(N)$\,,
let $Q_B$ be the quadratic form with matrix $B$ in the canonical
basis, and for all $v\in \MR^N$ let $D_v$ be the diagonal matrix with
diagonal $v$\,. For each $A\in \cala_N$ and every
$(\,u\vge (x \vge y))\in O^{N}_{\,1} \times O^{N}_{\,2}$\,, we thus have:
\[
\phi_A \,(u\vge (x\vge y)) 
\, = \,
Q_A(D_u\,x)+Q_A(D_u\,y)
\, = \,
Q_{B_{u\vge A}} (x) + Q_{B_{u\vge A}} (y) 
\]
where $B_{\,u\vge A}=D_u\,A\,D_u\in Sym(N)$\,, but also:
\[
\phi_A \,(u\vge (x\vge y)) 
\, = \,
Q_A(D_x\,u)+Q_A(D_y\,u)
\, = \,
Q_{B_{x\vge y \vge A}} (u)
\]
where  $B_{\,x\vge y \vge A}=B_{\,x\vge A}+B_{\,y\vge A}
=D_x\,A\,D_x+D_y\,A\,D_y$\,.
For all $B\in Sym(N)$\,, each critical point $u$ of 
$Q_B|_{O^{N}_{\,1}}$ satisfies:
$B\, u= \lambda\, u$ where $\lambda\in \MR$ 
is a Lagrange multiplier and we get: $Q_B(u)=\lambda$\,,
which proves that the maximum of $Q_B|_{O^{N}_{\,1}}$ is the largest 
eigenvalue of $B$\,. Similarly, let 
$Q_{B}^{[2]}\,(\,x\vge y)= Q_B(x)+Q_B(y)$ for all
$(\,x \vge y)\in  O^{N}_{\,2}$ and write the equations of the
submanifold $O^{N}_{\,2} \subset \MR^N \times \MR^N$ as:
$\langle\, x \vge x \rangle =1$\,,
$2\,\langle\, x \vge y \rangle =0$
and $\langle\, y \vge y \rangle =1$\,,
so that the critical points of $Q_{B}^{[2]}|_{O^{N}_{\,2}}$
are given by the Lagrange multipliers:
\[
\left \{ \begin{array}{l}
B\,x= a\,x+b\,y \\
B\,y=b\,x+c\,y
\end{array} \right.
\]
where $a\vge b \vge c \in \MR$\,, and
diagonalize the matrix
$\dsm
\left( \begin{array}{cc}
a &b \\
b &c
\end{array} \right)
$
in the orthonormal group to get an orthonormal basis  
$(\,\widetilde{x} \vge \widetilde{y}\,)$ of $\calp_{x\vge y}$ satisfying:
$Q_{B }^{[2]}\,(\,\widetilde{x}\vge \widetilde{y})
\,=\, Q_{B}^{  [2]}\,(\,x\vge y)$
and:
\[
\left \{ \begin{array}{l}
B\,\widetilde{x}= \alpha\,\widetilde{x} \\
B\,\widetilde{y}=\beta\,\widetilde{y}
\end{array} \right.
\]
where $\alpha \vge \beta \in \MR$ and
$Q_{B}^{[2]}\,(\,\widetilde{x}\vge \widetilde{y})=
\alpha+\beta$\,,
so the maximum of $Q_{B}^{ [2]}|_{O^{N}_{\,2}}$ is the sum of the two
largest eigenvalues of $B$\,.
Moreover, we get:
$A(\widetilde{x}\vge \widetilde{y})=A(x\vge y)$
and $B_{\,\widetilde{x}\vge \widetilde{y} \vge A}=B_{\,x\vge y\vge  A}$ 
since these matrices only depend on $\calp_{x\vge y}$\,,
so if $\Phi_{\,2}^N$ attains its maximum $\lambda_{\,2}^N$ at a point
$(u\vge (x_0 \vge y_0))\in O^{N}_{\,1} \times O^{N}_{\,2}$\,, we
get this way a point
$(u\vge (x \vge y))\in O^{N}_{\,1} \times O^{N}_{\,2}$ satisfying:
\[
\left \{ \begin{array}{ccl}
D_u\,A\,D_u\,x&= &\alpha\,x \\
D_u\,A\,D_u\,y&=&\beta\,y \\
\big(D_x\,A\,D_x+D_y\,A\,D_y)\,\,u&=&\lambda\,u
\end{array} \right.
\]
and $\lambda_{\,2}^N
=\phi_A(u\vge (x \vge y))=\alpha+\beta=\lambda$\,.
If $\alpha\,\beta \neq 0$ and if there exists an index
$1\leq i \leq N$ such that $u_i=0$\,, we get moreover:
$x_i=y_i=0$\,, and we saw that this implies 
$\lambda_{\,2}^N = \lambda_{\,2}^{N-1}$\,. 
Else, since we have $D_x\,u=D_u\,x$ and $D_y\,u=D_u\,y$
we get:
\[
\left \{ \begin{array}{ccl}
A\,D_x\,u&= &\alpha\,D_u^{-1}\,x \\
A\,D_y\,u&=&\beta\,D_u^{-1}\,y \\
D_x\,A\,D_x\,u+D_y\,A\,D_y\,u&=&(\alpha+\beta)\,u
\end{array} \right. \quad ,
\]
thus: 
$\alpha\,D_x\,D_u^{-1}\,x+\beta\,D_y\,D_u^{-1}\,y=
\alpha\,D_u^{-1\,}D_x\,x+\beta\,D_u^{-1}\,D_y\,y
=(\alpha+\beta)\,u$ since diagonal matrices commute,
and we obtain:
\[
(\alpha+\beta)\,u_i^2=\alpha\,x_i^2+\beta\, y_i^2
\quad \hbox{for all } 1\leq i \leq N \quad .
\]
Finally, for all $A\in \cala_N$\,, for all $u\in O^{N}_{\,1}$ and for all
$x\in O^{N}_{\,1}$ we get:
\[
Q_{B_{u\vge A}} (x) 
\, = \,
\sum_{1\,\leq\,\,i\vge j\,\,\leq \,N} a_{\, i\vge j}\,  
u_{\, i}\,  u_j\, x_i \, x_j
\, \leq \,
\Big(\sum_{1\,\leq\,i\,\leq \,N} \absol{u_i}\,\absol{x_i} \Big)^2
\, \leq \,
\norme{u}^2\,\norme{x}^2
\,=\,1
\]
by the Cauchy-Schwartz inequality, so all eigenvalues of $B_{u\vge A}$ 
are at most $1$\,.
But we have $\lambda_2^2=1$ and $(\lambda_{\,2}^N)_{\,N\in \MN}$ is
increasing, thus
if we assume $N\geq 3$ and $\alpha+\beta=\lambda_2^N > \lambda_2^{N-1}$
we obtain: $\alpha>0$ and $\beta>0$\,, hence  the above condition
$\alpha\,\beta \neq 0$ is fulfilled.

\section{The cases where $A$ is singular}\label{noninversible}

Assume that $N\geq 3$ and $\Phi_{\,2}^N$ attains its maximum 
$\lambda_2^N > \lambda_2^{N-1}$ at 
$(\,u\vge (x \vge y))\in O^{N}_{\,1} \times O^{N}_{\,2}$\,, thus we can
suppose that $u_i>0$ for all $1\leq i \leq N$\,, and let
$A=A(x\vge y)\in \cala_N$ be the matrix of signs.
If two rows (thus to columns) of $A$ are equal, we can assume by
the above symmetries that the last two lines of $A$ are equal to 
$(\,1 \,\, \cdots \,\, 1\,)$\,, and $M^{\ast}\in Sym(N-1)$ 
will denote thereafter the matrix obtained by removing the last line
and column of each symmetric matrix $M\in Sym(N)$ and we will set
$v^\ast =(\,v_{\,1}\,, \dots \,, v_{\,N-1}) \in \MR^{\,N-1}$ 
for each $v\in \MR^{\,N}$\,. We thus have
$A^{\ast \ast} \in \cala_{N-2}$ and:
\[
B_{\,  u\vge A}\,\,\,  = \,\,\,  
\left( \begin{array}{ccc|cc}
 &  & & u_{\,1}\,\,\,u_{N-1} & u_{\,1}\,\,\,u_{N} \\
\rule{0cm}{1cm} & B_{\,  u^{\ast\ast}\vge A^{\ast\ast}} 
& & \vdots & \vdots\\
\rule{0cm}{1cm}  &  & & u_{N-2}\,u_{ N-1} &u_{N-2}\,u_{ N}\\
&&&&\\
\hline
\rule{0cm}{0.6cm}u_{\,1}\,\,\,u_{N-1}  & 
\cdots & u_{N-2}\,u_{ N-1} & u_{N-1}^{\,\,2} & u_{N-1}\,u_{ N}\\
\rule{0cm}{0.6cm} u_{\,1}\,\,\,u_{N}   & \cdots & u_{N-2}\,\,\,\, u_{ N} 
&  u_{N-1}\,u_{ N} & u_{N}^{\,\,2}
 \end{array}\rule{0.2cm}{0cm}\right)
\] 
and since $u_N\neq 0$ we can set:
\[
V_u\,\,  = \,\,
\left( \begin{array}{ccc|cc}
 &  & & 0 &0\\
\rule{0cm}{0.8cm} & I_{N-2} & & \vdots & \vdots\\
\rule{0cm}{0.8cm}  &  & & 0 &0\\
\hline
\rule{0cm}{0.6cm} 0 & \cdots & 0 & \quad S_{\,u} & \\
0   & \cdots & 0 &  & 
 \end{array}\rule{0.1cm}{0cm}\right)
 \,\,\,  \hbox{ where }  \,\,\,  
S_{\,u}
\,\,  = \,\,  
\frac{1}{\sqrt{\, u_{N-1}^{\,\,2}+u_{N}^{\,\,2} \, }\,}
\,\,\,\left( \begin{array}{cc}
\rule{0cm}{0.5cm} u_{N-1} & u_{N}\\
\rule{0cm}{0.8cm} u_{N} & -\,  u_{N-1}\\
  \end{array}\rule{0.0cm}{0cm}\right) 
\] 
and obtain: \qquad   \qquad \qquad \qquad
$V_{\,u}^{\,\,-1}\,\,\,B_{u\vge A}\,\,\,  V_{\,u}\quad= $
\[
\left( \begin{array}{ccc|cc}
 &  & & u_{\,1}\,\,\,\sqrt{\, u_{N-1}^{\,\,2}+u_{N}^{\,\,2}} & 0 \\
\rule{0cm}{1cm} & B_{\,u^{\ast\ast}\vge \,A^{\ast\ast}} 
& & \vdots & \vdots\\
\rule{0cm}{1cm}  &  & 
& u_{N-2}\,\sqrt{\, u_{N-1}^{\,\,2}+u_{N}^{\,\,2}} & 0 \\
&&&&\\
\hline
\rule{0cm}{0.6cm}u_{\,1}\,\,\,\sqrt{\, u_{N-1}^{\,\,2}+u_{N}^{\,\,2}}  & 
\cdots & u_{N-2}\,\sqrt{\, u_{N-1}^{\,\,2}+u_{N}^{\,\,2}}
 & u_{N-1}^{\,\,2}+ u_{N}^{\,\,2} & 0\\
\rule{0cm}{0.6cm} 0   & \cdots &0
&  0 & 0
 \end{array}\rule{0.2cm}{0cm}\right)
\] 
\[
=\quad
\left( \begin{array}{c|c}
\rule{0cm}{0.6cm}
B_{\,\widetilde{u}^\ast\vge A^\ast} & 0_{\rule{0cm}{0.5cm}} \\
\hline
0 & 0
\end{array}\rule{0.0cm}{0cm}\right)
\quad \hbox{where } \,
\widetilde{u}^\ast 
\,=\,
\Big(\,u^{\ast\ast} \vge 
\sqrt{\, u_{N-1}^{\,\,2}+u_{N}^{\,\,2}\,}\,\,  \Big)
\, \in\, O^N_{\,1}
\quad .
\] 
The eigenvalues 
$b_1\geq \dots \geq b_{N}$ of 
$B_{\,u\vge A}$ are thus $0$ and the eigenvalues 
$d_1\geq \dots \geq d_{N-1}$
of $B_{\, \widetilde{u}^\ast\vge A^\ast}$ and there are two cases:
\smallskip

- either $d_2<0$ and we get: 
$\lambda_2^{N}=b_1+b_2=d_1+0\leq 1 \leq \lambda_2^{N-1}$\,,
\smallskip

- or $d_2\geq 0$ and we get:
$\lambda_2^{N}=b_1+b_2=d_1+d_2$\,,
but we also have:
\[
d_1+d_2 \,\,= \max_{(x^\ast \vge y^\ast)\, \in O^{N-1}_{\,2}} 
\phi_{A^\ast} (\,\widetilde{u}^\ast \vge (x^\ast \vge y^\ast)) 
\,\,\leq \max_{(x^\ast \vge y^\ast)\, \in O^{N-1}_{\,2}} 
\Phi_2^{N-1} (\,\widetilde{u}^\ast \vge (x^\ast \vge y^\ast)) 
\,\, \leq \lambda_2^{N-1}
\quad ,
\]
\noindent and in both cases we obtain
$\lambda_2^{N} = \lambda_2^{N-1}$ and we can conclude by induction.

\section{The case where $A$ is nonsingular}\label{inversible}

Assume that $N\geq 3$ and  $\Phi_N$ attains its maximum 
$\lambda_2^N > \lambda_2^{N-1}$ at 
$(\,u\vge (x^0 \vge y^0))\in O^{N}_{\,1} \times O^{N}_{\,2}$\,, 
so we can suppose that $u_k>0$ and $(x_k\vge y_k) \neq 0$ for all 
$1\leq k \leq N$\,. 
By changing the basis of $P_{x^0\vge y^0}$\,, then changing the sign of
each $(x_k\vge y_k)$\,, then permuting all the coordinates, we can obtain 
by the above symmetries 
$(x\vge y)\in O^{N}_{\,2}$ such that $\Phi_2^N(u\vge (x\vge y))=\lambda_2^N$
and:
\[
\qquad\,\,\,
x \,\,\, = \,\,\,
\big( \, R_0 \vge  R_s\,\sin \phi_s \vge \dots \vge R_1\,\sin \phi_1
\, \vge \, 
\,\, \rho_1\,\cos \psi_{1} \,\,  \vge  \, \dots \, \, \vge  
\,\,\rho_t\,\cos \psi_t \,\, \big)
\]
\[
\hbox{and} \quad
y \,\,\, = \,\,\,
\big( \,\, 0 \,  \vge  R_s\,\cos \phi_s \vge \dots \vge R_1\,\cos \phi_1
\, \vge \, 
- \rho_1\,\sin \psi_{1}  \vge  \, \dots \, \, \vge  
-\rho_t\,\sin \psi_t \, \big)
\]
where $s+t+1=N$\,, 
where we have:
$R_i>0$ and $0<\phi_i\leq \pi/2$ for each $1\leq i \leq s$ and 
$\rho_j>0$ and $0<\psi_j \leq \pi/2$ for each $1\leq j \leq t$\,,
and where the finite sequences $(\phi_i)_{1\leq i \leq s}$ and
$(\psi_j)_{1\leq j \leq t}$ are increasing.
Writing each $v\in \MR^N$ as:
\[
v \,\,\, = \,\,\,
\big( \, v^0 \vge  v^+_s \vge \dots \vge v^+_1
\, \vge \, 
v^-_1  \vge  \dots \vge  v^-_t \big)
\]
we get: \, $x^0\, x_k+y^0\, y_k=R_0\,x_k\geq 0$ \, 
for all $1 \leq k\leq N$
thus \,  $a_{1\vge k} (x\vge y)=+1$ \, ,
\[
x^+_i\, x^+_j+y^+_i\, y^+_j=R_i\,R_j\, \cos (\phi_i-\phi_j) > 0
\]
for all $1 \leq i\vge j \leq s$ since $-\pi/2 < \phi_i-\phi_j < \pi/2$
thus \, $a^{+\,+}_{i\vge j}=+1$ \, , \, and similarly:
\[
x^-_i\, x^-_j+y^-_i\, y^-_j=\rho_i\,\rho_j\, \cos (\psi_i-\psi_j) > 0
\]
for all $1 \leq i\vge j \leq t$ thus \, $a^{-\,-}_{i\vge j}=+1$ \, , \, 
and if $1 \leq i  \leq s$ and $1 \leq  j \leq t$ we get:
\[
x^+_i\, x^-_j+y^+_i\, y^-_j=R_i\,\rho_j\, \sin (\phi_i-\psi_j) 
\]
where \, $-\pi/2 < \phi_i-\psi_j < \pi/2$ \,
decreases with both indexes $k_i^+=s+2-i$ and $l_j^-=s+1+j$\,,
thus the symmetric sign matrix
$A(x\vge y)$ writes:
\[
A\,(x\vge y)\,\,\,= \,\,\,  
\left( \begin{array}{c|c|c}
 +1 & (+1) & (+1)  \\
\hline
\rule{0cm}{0.6cm}(+1) & (+1)  & C  \\
\hline
\rule{0cm}{0.6cm} (+1) & {}^t{} C & (+1)  \\
 \end{array}\right)
\] 
where $C\in M_{\,s\vge t}(\,\{-1\vge +1\}\,)$ ``\:has -1's bottom right'', 
that is if $c_{k_0\vge l_0}=-1$ we have: $c_{k\vge l}=-1$
for all $k\geq k_0$ and all $l\geq l_0$\,.
Since $A$ cannot have to equal lines due to the previous section, the
lines of $C$ must be distinct and distinct from $(1 \, \dots \, 1)$
thus $t\geq s$\,, as well as its columns 
thus $s\geq t$\,. We infer
$t=s$ hence $N=2\,s+1$\,, and (up to symmetries):
\[
A\,(x\vge y)\,\,\,=\,\,\,
A \,\,\,=\,\,\,
\left(\rule{0.0cm}{0cm}\begin{array}{c|cccc|cccc}
+\,1  & +\,1  & \cdots & \cdots&  +\,1& +\,1  & \cdots & \cdots &  +\,1 \\
\hline
\rule{0cm}{0.5cm} +\,1   &  +\,1  & \cdots  & \cdots  &  +\,1 & 
+\,1& \cdots & +\,1& -\,1\\
\rule{0cm}{0.5cm} \vdots  & \vdots & &  &\vdots
&\vdots &\iddots&\iddots& \vdots \\
\rule{0cm}{0.5cm} \vdots  & \vdots &&  &\vdots &
+\,1& \iddots &  & \vdots \\
\rule{0cm}{0.5cm} +\,1    & +\,1 & \cdots & \cdots & +\,1& 
  -\,1 & \cdots  & \cdots & -\,1  \\
\hline
\rule{0cm}{0.5cm}+\,1  &  +\,1   & \cdots &  +\,1 & -\,1 &  
+\,1  & \cdots  & \cdots  &  +\,1\\
\rule{0cm}{0.5cm}\vdots  & \vdots  & \iddots & \iddots & \vdots &  
\vdots &  &  &\vdots \\
\rule{0cm}{0.5cm}\vdots  &+\,1 &  \iddots &  & \vdots  & 
\vdots & &   &\vdots \\
\rule{0cm}{0.5cm}+\,1  &   -\,1  &   \cdots &  \cdots & -\,1 & 
+\,1  & \cdots  & \cdots  &  +\,1
\end{array}\rule{0.1cm}{0cm}\right) \quad,
\]
and since $a^{+\,-}_{i\vge j}$ is the sign of 
$x^+_i\, x^-_j+y^+_i\, y^-_j=R_i\,\rho_j\, \sin (\phi_i-\psi_j) 
$
where $-\pi/2 < \phi_i-\psi_j < \pi/2$ for all $1\leq i\vge j\leq s$
we can conclude that:
\[
0 \leq \phi_1 \leq \psi_1 \leq \cdots \leq \psi_{s-1} \leq 
\phi_s \leq \psi_s \leq\frac{\pi}{2} :=\phi_{s+1}
\quad .
\]
Moreover, if we had $\phi_i=\psi_i$ for a certain index $i$\,,
we could change the convention on the sign of 
$a^{+\,-}_{i\vge i}=0$ to obtain a matrix with to equal  lines which is
excluded, so we get finally:
\[
0 \leq \phi_1 < \psi_1 \leq \cdots < \psi_{s-1} \leq 
\phi_s < \psi_s \leq \phi_{s+1} =\frac{\pi}{2} 
\quad .
\]

\label{n=3}

If $s=1$ and $\lambda_2^3>\lambda_2^2=1$ we get thus:
\[
A\,\,\,=\,\,\,
\left(\rule{0.0cm}{0cm}\begin{array}{ccc}
+1  & +1  &  +1 \\
+1   &  +1  & -1\\
+1  &   -1   &  +1
\end{array}\rule{0.1cm}{0cm}\right) 
\quad \hbox{and} \quad 
B_{u\vge A}\,\,\,=\,\,\,
\left(\rule{0.0cm}{0cm}\begin{array}{ccc}
u_1^2  & u_1\,u_2  &  u_1\,u_3 \\
u_1\,u_2   &  u_2^2  & - u_2\,u_3\\
u_1\,u_3  &  - u_2\,u_3   &  u_3^2
\end{array}\rule{0.1cm}{0cm}\right) 
\]
for all $u\in O^N_{\,1}$, so the characteristic polynomial of 
$B_{u\vge A}$ is $P_u=X^3-X^2+4\,u_1^2\,u_2^2\,u_3^2$\,. 
Setting $P=X^3-X^2$ and
$\sigma=4\,u_1^2\,u_2^2\,u_3^2\in [\,0\vge 4/27\,]$
we get: $P_u=P+\sigma$ and the sum of all three roots of 
$P_u$ equals $1$\,, hence the sum $b_1+b_2$ of the largest two is maximal
when the least $b_3$ is minimal, that is when $\sigma=4/27$ is maximal, 
and we conclude that:
$\lambda_2^3=b_1+b_2=2/3+2/3=4/3$\,.

\section{The essential critical points}\label{essentiels}

If we have $N > 3$ and  $\Phi_2^N$ attains its maximum 
$\lambda_{\,2}^N > \lambda_{\,2}^{N-1}$ at 
$(\,u\vge (x \vge y))\in O^{N}_{\,1} \times O^{N}_{\,2}$\,, 
section \ref{inversible} proves that $N=2\,s+1$ is odd 
and gives the matrix $A=A(x\vge y)$ up to symmetries,
and in sections \ref{symetries} and \ref{extrema} 
we obtained after changing the basis of
$P_{x\vge y}$\,:
$D_u\,A\,D_u\,x=\alpha\,x$\,, $D_u\,A\,D_u\,y=\beta\,y$
and $\phi_A (u\vge (x\vge y ))=\alpha+\beta$\,, 
so we can assume that $1/3 < \beta \leq \alpha \leq 1$\,, and:  
$\alpha\, \big(x_k\big)^2+\beta\, \big(y_k\big)^2
=(\alpha+\beta)\,\big(u_k\big)^2$ where $u_k>0$
for all $1\leq k \leq N$\,, which
leads to:
\[
x \, = \,\,
\sqrt{\frac{\alpha+\beta}{\alpha}} \,\,
\Big( \, u^0 \,\sin \theta_{s+1}\vge  u^+_s\,\sin \theta_s \vge 
\dots \vge u_1^+\,\sin \theta_1
\, \vge \, 
\,\, u^-_1\,\cos \varphi_1  \,  \vge  \, \dots \, \, \vge  
\,\, u^-_{s}\,\cos \varphi_s 
\, \Big)
\]
\[
y \,= \,\,
\sqrt{\frac{\alpha+\beta}{\beta}} \,\,
\Big( \,  u^0 \,\cos \theta_{s+1}  \,\vge \, u^+_s \,\cos\theta_s
\vge \dots \vge  u_1^+\,\cos\theta_1
 \vge  -u^-_1\,\sin\varphi_1   \vge  \dots  \vge  
-u^-_{s}\,\sin\varphi_s \, \Big)
\]
where after this change of basis of $P_{x\vge y}$
(which writes $\phi_k \mapsto
\phi_k+\phi_0$ and $\psi_k \mapsto \psi_k-\phi_0$
for each $1\leq k \leq s$)  
and this affinity in each plane $(y_k\vge x_k)$
we get:
\[
\theta_1 < \theta_2 < \cdots  < \theta_s  < \theta_{s+1}  
< \theta_1+\pi 
\]
since half-planes are preserved by linear transformations. 
We also have:
\[
A^{\,-1}\,\,\,=\,\,\,\frac{1}{2}\,\,\,
\left(\rule{0.0cm}{0cm}\begin{array}{c|cccc|cccc}
0  & 0  & \cdots & 0 &  +\,1&  0 & \cdots &    0 & +\,1 \\
\hline
\rule{0cm}{0.5cm} 0   & 0  & \cdots  & \cdots  &  0 & 
(0)&  & +\,1& -\,1\\
\rule{0cm}{0.5cm} \vdots  & \vdots & &  &\vdots
&  &\iddots&\iddots&  \\
\rule{0cm}{0.5cm} 0  & \vdots &&  &\vdots &
+\,1& \iddots &  (0)  &  \\
\rule{0cm}{0.5cm} +\,1    & 0 & \cdots & \cdots & 0& 
  -\,1 &   &  & \\
\hline
\rule{0cm}{0.5cm}0  &  (0)   &  &  +\,1 & -\,1 &  
0  & \cdots  & \cdots  &  0 \\
\rule{0cm}{0.5cm}\vdots  &   & \iddots & \iddots &  &  
\vdots &  &  &\vdots \\
\rule{0cm}{0.5cm}0  &+\,1 &  \iddots & (0) &   & 
\vdots & &   &\vdots \\
\rule{0cm}{0.5cm}+\,1  &   -\,1  &    &   &  & 
0  & \cdots  & \cdots  &  0
\end{array}\rule{0.1cm}{0cm}\right)
\]
and the equation 
\, $D_u\,A\,D_u\,v=\lambda\, v$ \, also writes: \quad 
$D_u\,v=\lambda\, A^{\,-1}\,D_u^{\,-1}\, v$ \, , \,
thus if we set: 
  \,  $\dsm \varphi_0 = \theta_{s+1} - \frac{\pi}{2}$ \,
, \,  $u_0^-=u^0$ \, 
, \,  $\dsm \theta_0 = \varphi_{s} - \frac{\pi}{2}$ \,
, \,  $\dsm \varphi_{-1} = \theta_{s} - \frac{\pi}{2}$ \,
and \, $u_0^+=u_s^-$ \,,
we infer from 
$D_u\,A\,D_u\,x=\alpha\,x$ and $D_u\,A\,D_u\,y=\beta\,y$
 the two systems of equations for $0\leq k \leq s$\,:
\[
\left \{ \begin{array}{l}
2\,\big(u^+_k\big)^2\, \sin \theta_k
=
\alpha\,(\cos \varphi_{k-1}-\cos \varphi_{k})  \\
2\,\big(u^+_k\big)^2\, \cos \theta_k
=
\beta\,(\sin \varphi_{k}-\sin \varphi_{k-1}) 
\end{array} \right . \qquad 
\]
\[
\left \{ \begin{array}{l}
2\,\big(u^-_k\big)^2\, \cos \varphi_k
=
\alpha\,(\sin \theta_{k+1}-\sin \theta_k)\\
2\,\big(u^-_k\big)^2\, \sin \varphi_k
=
\beta\,(\cos \theta_k-\cos \theta_{k+1} )
\end{array} \right . 
\quad . \, \,
\]
For each $0\leq k \leq s$\,, the first system implies:
\[
\alpha\, \cos \theta_k \,(\cos \varphi_{k-1}-\cos \varphi_{k}) 
(\sin \varphi_{k}+\sin \varphi_{k-1}) 
=
\beta \, \sin \theta_k \,(\sin \varphi_{k}-\sin \varphi_{k-1}) 
(\sin \varphi_{k}+\sin \varphi_{k-1}) 
\]
hence (since $\sin^2 \varphi_{k}-\sin^2 \varphi_{k-1}
\,=\,\cos^2 \varphi_{k-1}-\cos^2 \varphi_{k}$):
\[
\alpha\, \cos \theta_k \,(\cos \varphi_{k-1}-\cos \varphi_{k}) 
(\sin \varphi_{k}+\sin \varphi_{k-1}) 
=
\beta \, \sin \theta_k \,(\cos \varphi_{k-1}-\cos \varphi_{k}) 
(\cos \varphi_{k-1}+\cos \varphi_{k}) 
\]
thus either \, $\cos \varphi_{k-1}=\cos \varphi_{k}$ \, or:
\[
\alpha\, \cos \theta_k \, (\sin \varphi_{k}+\sin \varphi_{k-1}) 
=
\beta \, \sin \theta_k \, (\cos \varphi_{k-1}+\cos \varphi_{k}) 
\quad ,
\]
but in the first case, we get $\sin \theta_k=0$
thus $\cos \theta_k \neq 0$\,, hence 
\, $\sin \varphi_{k-1}\neq\sin \varphi_{k}$ \,
thus \, $\sin \varphi_{k-1}=-\sin \varphi_{k}$ \, , \,
hence the second equation is always fulfilled, and it also writes:
\[
\beta\, \sin \theta_{k} \, \cos \varphi_{k-1} -
\alpha\,  \cos \theta_{k} \, \sin \varphi_{k-1}
\,=\,
-\beta\, \sin \theta_{k} \, \cos \varphi_{k} +
\alpha\,  \cos \theta_{k} \, \sin \varphi_{k}
\quad .
\]
Similarly, the second system implies 
for each $0\leq k \leq s$\,:
\[
\beta\, \sin \theta_{k+1} \, \cos \varphi_{k} -
\alpha\,  \cos \theta_{k+1} \, \sin \varphi_{k}
\,=\,
-\beta\, \sin \theta_{k} \, \cos \varphi_{k} +
\alpha\,  \cos \theta_{k} \, \sin \varphi_{k}
\quad ,
\]
and we infer from these equations: \qquad 
$-\beta\, \sin \theta_{k} \, \cos \varphi_{k} +
\alpha\,  \cos \theta_{k} \, \sin \varphi_{k}
\quad = $
\[
\beta\, \sin \theta_{k+1} \, \cos \varphi_{k} -
\alpha\,  \cos \theta_{k+1} \, \sin \varphi_{k}
\,=\,
\beta\, \sin \theta_{k} \, \cos \varphi_{k-1} -
\alpha\,  \cos \theta_{k} \, \sin \varphi_{k-1}
\,=:\, h \in \MR \,.
\]
But  for each $0\leq k \leq s$ we obtain by the first system:
\[
\left \{ \begin{array}{l}
2\,\big(u^+_k\big)^2\, \sin^2 \theta_k
=
\alpha\,(\cos \varphi_{k-1}-\cos \varphi_{k})\, \sin \theta_k  \\
2\,\big(u^+_k\big)^2\, \cos^2 \theta_k
=
\beta\,(\sin \varphi_{k}-\sin \varphi_{k-1})\, \cos \theta_k 
\end{array} \right . \qquad
\]
hence:  \qquad
$\dsm
2\,\big(u_k^+\big)^2\, 
\Big(\,\frac{\beta}{\alpha} \,\sin^2 \theta_{k} 
+\frac{\alpha}{\beta} \, \cos^2 \theta_{k} \, \Big)
\, = \, 2\,h
$ \quad , \quad and similarly by the second:
\[
2\,\big(u^-_k\big)^2\, 
\Big( \, \frac{\beta}{\alpha} \,\cos^2 \varphi_{k} 
+\frac{\alpha}{\beta} \, \sin^2 \varphi_{k} \, \Big) 
\, = \, 2\,h
\quad , \qquad \qquad \qquad \qquad \qquad \, \, \,
\]
thus if $0\leq k \leq s$ we have: \qquad
$\dsm
\frac{\beta\,\big(x_k^+\big)^2}{\alpha+\beta}
+\frac{\alpha\,\big(y_k^+\big)^2}{\alpha+\beta}
\, = \, 
h
$ \qquad
and if $1\leq k \leq s$:
\[
\frac{\beta\,\big(x_k^-\big)^2}{\alpha+\beta}
+\frac{\alpha\,\big(y_k^-\big)^2}{\alpha+\beta}
\, = \, 
h \qquad ,
\]
thus it comes by adding these equations: \,
$\beta\, \norme{x}^2+\alpha\,\norme{y}^2 \, = \,
(2\,s+1)\,(\alpha+\beta)\,h$\,,  but we must 
have: \,
$\norme{x}\,=\,\norme{y}\,=\,1$ \, thus we finally infer:
\, $\dsm h\,=\,\frac{1}{2\,s+1} $ \, .
\smallskip

\noindent
But we also have: \quad
$
\dsm \norme{x}^2
\,=\,
\frac{\alpha+\beta}{\alpha} \, \Big( \, 
\sum_{k=1}^s \big(u^+_k\big)^2\, \sin^2 \theta_k  
\, +\,
\sum_{k=0}^s \big(u^-_k\big)^2\, \cos^2 \varphi_k 
\, \Big) 
$ \, , \, thus by the systems above we infer:
\[
\frac{2}{\alpha+\beta}\,\,\norme{x}^2
\,=\,
\sum_{k=1}^s 
\big(\,\cos \varphi_{k-1}-\cos \varphi_{k}\,\big)\, \sin \theta_k
\, +\,
\sum_{k=0}^s 
\big(\,\sin \theta_{k+1}-\sin \theta_k\,\big)\, \cos \varphi_k
\quad ,
\]
hence $\norme{x}=1$ implies (since $\cos \varphi_0=\sin \theta_{s+1}$
and $\sin \theta_0=-\cos \varphi_s$):
\[
\sin \theta_{1} \, \sin \theta_{s+1}
\,+\,
\sum_{k=1}^{s} 
\big(\,\sin \theta_{k+1}-\sin \theta_k\,\big)\, \cos \varphi_k
\,=\,
\frac{1}{\alpha+\beta}
\quad .
\]

Moreover, we obtained for each $1\leq k \leq s$\,:
\[
\left\{\begin{array}{lcc}
-\cos \theta_{k+1} \,\big(\,\alpha\, \sin \varphi_k\,\big)
\,+\,\sin \theta_{k+1} \,\big(\, \beta\, \cos \varphi_k\,\big)
&=&\dsm\frac{1}{2\,s+1}\\
\qquad\cos \theta_{k} \,\,\big(\,\alpha\, \sin \varphi_k\,\big)
\,\,-\,\,\sin \theta_{k} \,\,\,\big(\, \beta\, \cos \varphi_k\,\big)
&=&\rule{0cm}{0.7cm} \dsm\frac{1}{2\,s+1}
\end{array}\right.
\]
and this system in 
$(\alpha\,\sin\varphi_k \vge \beta\, \cos \varphi_k)$
has determinant $\sin\big(\theta_k-\theta_{k+1}\big)<0$
so it comes:
\[
\left(\begin{array}{c}
\alpha\,\sin \varphi_k \\
\beta\,\cos \varphi_k 
\end{array}\right)
=
\frac{1}{(2s+1)\, \sin(\theta_{k+1}-\theta_k)}\,
\left(\begin{array}{cc}
\sin \theta_k+\sin \theta_{k+1}\\
\cos \theta_k+\cos\theta_{k+1}
\end{array}\right) \quad ,
\]
which implies (since $\cos^2 \varphi_k + \sin^2 \varphi_k=1$)\,:
\[
\alpha^2\,(\cos\theta_{k+1}+\cos \theta_k)^2 +
\beta^2\,(\sin \theta_{k+1}+\sin \theta_k)^2=
(2s+1)^2\,\alpha^2\,\beta^2 \sin^2(\theta_{k+1}-\theta_k)\quad ,
\]
hence by the usual trigonometric formulae:
\[
\alpha^2\,\cos^2\Big(\,\frac{\theta_{k}+\theta_{k+1}}{2}\,\Big) 
+\beta^2\,\sin^2\Big(\,\frac{\theta_{k}+\theta_{k+1}}{2}\,\Big) 
=
(2s+1)^2\,\alpha^2\,\beta^2 \,
\sin^2\Big(\,\frac{\theta_{k}-\theta_{k+1}}{2}\,\Big) 
\]
since
$\dsm 0 < \frac{\theta_{k+1}-\theta_{k}}{2} < \frac{\pi}{2}$\,, 
hence
$\dsm \cos \Big( \frac{\theta_{k}-\theta_{k+1}}{2}\Big)\neq 0$\,, 
and we get the induction relation:
\[
\sin \Big(\,\frac{\theta_{k+1}-\theta_{k}}{2}\,\Big)
=\frac{1}{2\,s+1}\,\sqrt{\,
\frac{1}{\beta^2}
-\Big(\,\frac{1}{\beta^2}-\frac{1}{\alpha^2}\,\Big)\,
\sin^2\Big(\,\frac{\theta_{k}+\theta_{k+1}}{2}\,\Big)\,}
\]
which determines $\theta_{k+1}$ knowing  $\theta_k$
(see section \ref{equadiff}).
Moreover, we obtained the relation:
\[
\beta\,   \sin \theta_{1}\, \sin \theta_{s+1}  +
\alpha\,   \cos \theta_1  \,  \cos \theta_{s+1} 
=\frac{1}{2s+1}
\]
which is a boundary condition on the sequence
$(\theta_k)_{\,1\leq k \leq s+1}$\,, and the ``norm equation'': 
\[
\sin \theta_{1}\, \sin \theta_{s+1}  +
\sum_{k=1}^{s}\,\,
\big(\sin \theta_{k+1}-\sin \theta_k\big)\, \cos \varphi_k
\,=\,
\frac{1}{\alpha+\beta}
\quad ,
\]
hence by the above expression of $\beta\,\cos \varphi_k$\,:
\[
\sin \theta_{1}\, \sin \theta_{s+1}  +
\frac{2}{(2s+1)\,\beta}
\sum_{k=1}^{s}\,\,
\cos^2\Big(\frac{\theta_k+\theta_{k+1}}{2}\Big) 
\,=\,
\frac{1}{\alpha+\beta} \quad ,
\]
and, by usual trigonometric formulae, 
we can finally write the norm equation as:
\[
\frac{1}{2\,s+1}\,\,\sum_{k=1}^{s} \,\, 
\cos\big(\,\theta_k+\theta_{k+1}\,\big) 
\,=\,
\frac{(s+1)\,\beta-s\,\alpha}{(2s+1)\,(\alpha+\beta)} 
-\beta\,\sin \theta_{1}\, \sin \theta_{s+1}  
\quad .
\]
Moreover, these three relations are invariant under the global
translation $\theta_k \mapsto \theta_k+\pi$ for $1\leq k \leq s+1$\,, 
thus we can assume hereafter
that: \, $0 \leq \theta_1 <\pi$ \, .

\section{The differential equation}\label{equadiff}

If we have $s \geq 2$ and  
$\lambda_{\,2}^{2\,s+1} > \lambda_{\,2}^{2\,s-1}$\,,
we obtained two real numbers \,
$1/3 < \beta \leq \alpha \leq 1$ satisfying \,
$\alpha + \beta > 4/3$ \, and a  finite sequence:
\[
0 \leq \theta_1 < \theta_2 < \cdots  < \theta_s  < \theta_{s+1}  
< \theta_1+\pi < 2\,\pi
\]
satisfying the induction relation for each $1\leq k \leq s$\,:
\[
\sin \Big(\,\frac{\theta_{k+1}-\theta_{k}}{2}\,\Big)
\,=\,
\frac{1}{2\,s+1}\,\sqrt{\,\frac{1}{\beta^2}
-\Big(\,\frac{1}{\beta^2}-\frac{1}{\alpha^2}\,\Big)\,
\sin^2\Big(\,\frac{\theta_{k}+\theta_{k+1}}{2}\,\Big)\,}
\quad ,
\]
hence if $s$ is large enough, each $\theta_{k+1}-\theta_{k}$ is small 
and it becomes:
\[
\theta_{k+1}-\theta_{k}
\, \simeq \,
\frac{2}{2\,s+1}\,\sqrt{\,
\frac{1}{\beta^2}
-\Big(\,\frac{1}{\beta^2}-\frac{1}{\alpha^2}\,\Big)\,
\sin^2\Big(\,\frac{\theta_{k}+\theta_{k+1}}{2}\,\Big)\,}
\quad ,
\]
thus each $\theta_k$ approximates 
$\dsm y\Big(\,\frac{2\,k-2}{2\,s+1}\,\Big)$
by some kind of ``middle-point at the goal method'',
where  $y\dpe [\,0\vge 1\,] \dans \MR$ satisfies the
the differential equation: 
\[
y' \,=\, f(y) \,:=\,
\frac{1}{\beta}\,\sqrt{\,1
-\Big(1-\frac{\beta^2}{\alpha^2}\,\Big)\,\sin^2 y\,} 
\]
and the initial condition $y(0)=\theta_1$\,.
First of all, we compute easily: 
\[
\norme{f}_\infty  = \, \frac{1}{\beta}  \,\leq  \,3
\quad ,\quad  
\norme{f'\,}_\infty 
 = \,
\frac{1}{\beta} - \frac{1}{\alpha}
\, \leq \,
2
\quad \hbox{and} \quad 
\norme{f''\,}_\infty 
\, = \,
\alpha\,\Big(\frac{1}{\beta^2} - \frac{1}{\alpha^2}\Big)
\, \leq \,
8
\quad , 
\]
thus $f$ is $2$-Lipschitz.
For each $\theta\in \MR$\,, let:  
\,
$ g_\theta (x) \, = \, \sin x - h\,f (\theta+x )$
\,
for all $x\in \MR$ where $\dsm h=1/(2\,s+1)$\,,  
thus we have: \,
$
g_\theta'(x) \, = \, \cos x - h\,f' (\theta+x )
\, \geq \,\cos x - 2\,h
$ \, , \,
hence 
$g_\theta$ is strictly increasing on $[\,0\vge x_0\,]$ where 
$ x_0 \, = \, \arccos ( 2\,h )$\,, 
and we have: \,  $g_\theta(0)\,<\,0$ \,  and  
$ g_\theta(x)\, \geq \,  \sqrt{1-4\,h^2}  -3\,h >0$ \,  
if $s\geq 2$ and 
$ x\in [ \, x_0\vge \frac{\pi}{2} \, ]$\,, 
thus  there exists a unique point
$\delta \theta\in [\,0\vge \frac{\pi}{2}\,]$ such that: \, 
$g_\theta(\delta \theta)=0$ \, . \, 
This shows that the induction relation
determines  $\theta_{k+1}$ knowing  $\theta_k$\,, and we get: \,
$\theta_{k+1}=\,\theta_k + 2\,\delta \theta_k$ \, 
where \, $0<\delta \theta_k \leq \arcsin (3\,h)$ \, .
\bigskip

For each $\varphi_0 \in \MR$\,, consider on one hand:
\, $\widetilde{\varphi_1}=\varphi_0+\delta \varphi_0$ \,
and
\, $\widetilde{\varphi_2}=\varphi_0+2\,\delta \varphi_0$ \, 
as above,
and on the other hand,
let $y$ be the solution of \, $y'=f(y)$ \, satisfying 
\, $y(0)=\varphi_0$ \, 
and set: \, $\varphi_1= y(h)$ \, and \, $\varphi_2= y(2\,h)$  \, . \, 
We obtain:
\[
y''
\, = \,  
\frac{1}{2}\,\Big(\frac{1}{\alpha^2}-\frac{1}{\beta^2}\Big)\,
\sin\big(2\,y\big)
\quad , \quad
y'''
\, = \,  
\Big(\frac{1}{\alpha^2}-\frac{1}{\beta^2}\Big)\,
\cos\big(2\,y\big)\, y'
\qquad \hbox{and}
\]
\[
y^{(4)}
\, = \,  
\frac{1}{4}\,\Big(\frac{1}{\alpha^2}-\frac{1}{\beta^2}\Big)^2\,
\sin\big(4\,y\big)
-2\, \Big(\frac{1}{\alpha^2}-\frac{1}{\beta^2}\Big)\,
\sin\big(2\,y\big)\, \big(y'\big)^2
\quad , 
\]
thus: \, 
$\norme{y'}_\infty \,  \leq\, 3$ \, , \,
$\norme{y''}_\infty \, \leq \, 4$ \, , \,
$\norme{y'''}_\infty \, \leq\,24$ \, and \, 
$\norme{y^{(4)}}_\infty \,  \leq\, 160$ \, . \, 
We infer at first that: \,
$\varphi_0 \,\leq\, \varphi_2 \,\leq\, 6\,h$ \, , \, but if we write:
\[
\varphi_0=\varphi_1-h\,f(\varphi_1)+\frac{h^2}{2}\,y'' (h)
-\frac{h^3}{6}\,y''' (t_1)
\, \hbox{ and } \,
\varphi_2=\varphi_1+h\,f(\varphi_1)+\frac{h^2}{2}\,y'' (h)
+\frac{h^3}{6}\,y''' (t_2)
\]
where \, $0< t_1 < h < t_2 < 2\,h$ \, , \, we obtain moreover:
\, $\dsm \absol{\varphi_{2}-\varphi_{0}-2\,h\, f(\varphi_{1})}
\, \leq \, 8 \, h^3
$ \, , \,
and we have: \,  $f(\varphi_{1}) \geq 1/\alpha \geq 1$ \, , \,
thus if we set:
\[
\Delta \varphi_0 \,=\, \frac{\varphi_2-\varphi_0}{2}
\hbox{ , \, we get:} \quad
0.9958\,h\, f(\varphi_{1})  \,\leq\, \Delta \varphi_{0} \,\leq\, 1.0042\,h\, f(\varphi_{1}) 
\, \hbox{ whenever } \, s\geq 15
\]
that is the essence of the quadratic convergence of the middle-point method,
and from now on we will assume that \,  $s\geq 15$ \, . 
In the same way, we have first:
\[
\varphi_0=\varphi_1-h\,f(\varphi_1)+\frac{h^2}{2}\,y'' (t_3)
\quad
\hbox{and}
\quad
\varphi_2=\varphi_1+h\,f(\varphi_1)+\frac{h^2}{2}\,y'' (t_4)
\]
where \, $0< t_3 < h < t_4 < 2\,h$ \, , \, hence: \,
$\dsm \Absole{\frac{\varphi_0+\varphi_2}{2}-\varphi_1}
\,\leq\, 2\,h^2
$ \, , \, but if we write:
\[
\qquad \,\,\, 
\varphi_0=\varphi_1-h\,f(\varphi_1)+\frac{h^2}{2}\,y'' (h)
-\frac{h^3}{6}\,y''' (h)+\frac{h^4}{24}\,y^{(4)} (t_5)
\]
\[
\hbox{and} \quad
\varphi_2=\varphi_1+h\,f(\varphi_1)+\frac{h^2}{2}\,y'' (h)
+\frac{h^3}{6}\,y''' (h)+\frac{h^4}{24}\,y^{(4)} (t_6) 
\]
where \, $0< t_5 < h < t_6 < 2\,h$ \, , \, we obtain moreover:
\[
\Absole{\frac{\varphi_0+\varphi_2}{2}-\varphi_1
-\frac{1}{4}\,\Big(\frac{1}{\alpha^2}-\frac{1}{\beta^2}\Big)\,
\sin\big(2\,\varphi_1\big)\,h^2}
\,\leq\,
\frac{20}{3}\,h^4
\quad ,
\]
thus since $f$ is $2$-Lipschitz:
\[
\Absole{f \Big(\,\frac{\varphi_0+\varphi_2}{2}\, \Big)
-f \Big(\,\varphi_1
-\frac{1}{4}\,\Big(\frac{1}{\alpha^2}-\frac{1}{\beta^2}\Big)\,
\sin\big(2\,\varphi_1\big)\,h^2\, \Big)}
\,\leq\,
\frac{40}{3}\,h^4
\quad ,
\]
hence since \, $\norme{f''}_\infty \,\leq 8$ \, :
\[
\Absole{f \Big(\,\frac{\varphi_0+\varphi_2}{2}\, \Big)
-f (\varphi_1)
+\frac{1}{4}\,\Big(\frac{1}{\alpha^2}-\frac{1}{\beta^2}\Big)\,
f' (\varphi_1)\,\sin\big(2\,\varphi_1\big)\,h^2}
\,\leq\,
\frac{88}{3}\,h^4 
\quad .
\]
Similarly, we obtained above:
\, $\dsm \absol{\Delta \varphi_{0}-h\, f(\varphi_{1})}
\, \leq \, 4 \, h^3
$ \, , \,
hence since $f$ is $2$-Lipschitz:
\[
\Absole{ \Delta \varphi_0
-h\,f \Big(\,\frac{\varphi_0+\varphi_2}{2}\, \Big)}
\, \leq\,
4\,h^3+
h\,\Absole{ f \Big(\,\frac{\varphi_0+\varphi_2}{2}\, \Big)
-f(\varphi_{1})}
\, \leq\,
8\,h^3
\quad ,
\]
as well as:
\, $0\,\leq\,\Delta \varphi_0 \,\leq\,  3\,h$ \, 
hence: 
\, $\dsm
\Absole{\sin \Delta \varphi_{0}-\Delta \varphi_{0}}
\, \leq \,
\frac{9}{2}\,h^3
$ \, and we infer at first:
\[
\Absole{g_{\varphi_0}\big(\,\Delta \varphi_0\,\big)}
\, \leq\, 
\frac{25}{2}\,h^3
\quad ,
\]
but we get now:
\, $\dsm
\Absole{ \frac{\varphi_2-\varphi_0}{2}
-f(\varphi_1)\,h-\frac{1}{6}\,
\Big(\frac{1}{\alpha^2}-\frac{1}{\beta^2}\Big)\,
\cos\big(2\,\varphi_1\big)\, f(\varphi_1)\, h^3 }
\, \leq\, \frac{20}{3}\,h^4
$ \,
thus:
\[
\Absole{ \Delta \varphi_0
-h\,f \Big(\,\frac{\varphi_0+\varphi_2}{2}\, \Big)}
\, \leq\,
\Big(\frac{1}{\beta^2}-\frac{1}{\alpha^2}\Big)\,
\Absole{\frac{\cos(2\,\varphi_1)\, f(\varphi_1)}{6}
+\frac{\sin(2\,\varphi_1)\,f' (\varphi_1)}{4}}\,h^3
+ 36\,h^4
\quad ,
\]
but we have: 
\, $\dsm
\Absole{\sin \Delta \varphi_{0}-\Delta \varphi_{0}}
\, \leq \,
\frac{1}{6}\,\big(\,1.0042\,h\, f(\varphi_{1})\,\big)^3 
\, \leq \,
0.17\, f^3(\varphi_{1})
$ \, and we obtain:
\[
\Absole{g_{\varphi_0}\big(\,\Delta \varphi_0\,\big)}
 \leq
\Big(\, 
\Big(\frac{1}{\beta^2}-\frac{1}{\alpha^2}\Big)\,
\Absole{\frac{\cos(2\varphi_1)\, f(\varphi_1)}{6}
+\frac{\sin(2\varphi_1)\,f' (\varphi_1)}{4}}
+0.17 f^3(\varphi_{1})+36h \Big)\,h^3
\]
hence:
$ \,
\absol{g_{\varphi_0}\big(\,\Delta \varphi_0\,\big)}
\, \leq
9.0553\,h^3
$ \,.
Finally, we have to estimate the derivative of $g_{\varphi_0}$\,: we have:
\, $0 \,\leq\, \delta \varphi_{0} \,\leq\, \arcsin(3\,h)$ \, 
and:
\, $0\,\leq\,\Delta \varphi_0 \,\leq\, 3\,h \,\leq\, \arcsin(3\,h)$ \, 
, \, hence for all real
$x\in [\,\delta\varphi_0 \vge \Delta\varphi_0\,]$ we obtain:
\, $ \absol{g_{\varphi_0} ' (x)}
\,\geq\, \sqrt{1-9\,h^2}- 2\,h
\,\geq\, 0.9307 $ \, , \,
and since we have:
\, $g_{\varphi_0} (\delta\varphi_0)=0$ \, 
we get:
\[
\Absole{\Delta \varphi_{0}-\delta\varphi_0 }
\, \leq \,
\frac{9.0553}{0.9307}\,\,h^3
\, \leq \,
9.73\,h^3 \quad .
\]
We can conclude that:
$
\absol{\varphi_2-\widetilde{\varphi_2}}
\, \leq\, 
19.46\, h^3
$ 
and the "consistency error in the middle" is:
\[
\Absole{\varphi_1-\frac{\varphi_0+\widetilde{\varphi_2}}{2}}
\, \leq\, 
\Absole{\varphi_1-\frac{\varphi_0+\varphi_2}{2}}
+
\Absole{\Delta \varphi_{0}-\delta\varphi_0 }
\, \leq \,
2\,h^2+9.73\,h^3 
\,\leq\, 2.314\,h^2  
\quad .
\]

For all integers $0\leq j \leq 2\,s+1$ let \, $t_j=j\,h$\,, \,and for all
integers $1\leq k \leq s$ let $y_k$ be the solution of $y'=f(y)$ 
such that \, $y_k(t_{2\,k-2})=\theta_k$ \, and let
\, $\widetilde{\theta}_{k+1}=y_k(t_{2\,k})$ \, as above.
Set \, $y=y_1$ \, , \,  and for all \, $1\leq k \leq s+1$ \, set:
\, $
\eps_k \,=\, \absol{y(t_{2\,k-2})-\theta_k} \,\,=\,\,
\absol{(y-y_k)\,(t_{2\,k-2})}
$ \,
to obtain: \, $\eps_1=0$ \, and \, 
$\absol{(y-y_k)\,'}\, \leq 2\,\absol{y-y_k}$ \, , \, thus
by the Gronwall lemma:
\[
\absol{y(t_{2\,k})-\widetilde{\theta}_{k+1}}
\,=\,
\absol{(y-y_k)\,(t_{2\,k})}\,
\,\leq \,
\eps_k\,e^{\,2\,\absol{\,t_{2\,k}-t_{2\,k-2}\,}}
\,=\,
e^{4\,h}\,\eps_k
\]
for all integers $1\leq k \leq s$\,, hence it comes:
\[
\eps_{k+1}
\,\leq \,
\absol{y(t_{2\,k})-\widetilde{\theta}_{k+1}}
+
\absol{\widetilde{\theta}_{k+1}-\theta_{k+1}}
\,\leq \,
e^{4\,h}\,\eps_k+ 19.46\,h^3 
\]
by the above estimate, thus for all integers $1\leq k \leq s+1$
we get:
\[
\Absole{\theta_k-y(t_{2\,k-2})} 
\,\,=\,\,
\eps_{k}
\, \leq \,
19.46\,h^3 \, \frac{e^{\,4\,h\,(k-1)}-1}{e^{\,4\,h}-1}
\, \leq \,
4.865\,h^2 \, \big( e^{\,2\,t_{2\,k-2}}-1\big)  
\quad ,
\]
where $y$ is the solution of \, $y'=f(y)$ \, satisfying 
\, $y(0)=\theta_1$ \, . \, 
Moreover, we also get:
\, $\dsm
\absol{(y-y_k)\,(t_{2\,k-1})}\,
\,\leq \,
e^{\,2\,h}\,\eps_k
$ \,
for all $1\leq k \leq s$ \,  and we obtained:
\[
\Absole{y_k (t_{2\,k-1})-\frac{\theta_k+\theta_{k+1}}{2}}
\, \leq \,
2\,h^2+9.73\,h^3 
\quad ,
\]
thus it comes:
\[
\Absole{\frac{\theta_k+\theta_{k+1}}{2}-y (t_{2\,k-1})}
\, \leq \,
2\,h^2+9.73\,h^3  + e^{\,2\,h}\,\eps_k
\, \leq \,
2\,h^2+9.73\,h^3  +
5\,h^2 \, \big( e^{(4\,k-2)\,h}-1-2\,h\big)
\]
hence: \,
$ \dsm
\Absole{\frac{\theta_k+\theta_{k+1}}{2}-y (t_{2\,k-1})}
\, \leq \,
5\,h^2 \, \big( e^{\,2\,t_{2\,k-1}}-1\big)+2\,h^2 
$ \, .

\section{The boundary condition}\label{cb}

When $\theta_1$ is fixed, 
$\theta_{s+1}$ goes to $y(1)$ when $s$ goes to infinity, 
and the boundary condition:
\[
g(\theta_{s+1})\,:=\,\alpha\,   \cos \theta_1  \,  \cos \theta_{s+1} 
+ \beta\,   \sin \theta_{1}\, \sin \theta_{s+1} 
\,=\,\frac{1}{2\,s+1}\,=\,h 
\]
where $\theta_1<\theta_{s+1}<\theta_1+\pi$
becomes: \, $y(1)=\theta_\perp$ \, where
$g(\theta_\perp)=0$ \,
and \, $\theta_1<\theta_\perp\leq \theta_1+\pi$\,, thus if we set:
\, $\theta_\perp=\theta_1+\pi/2+\phi$  \,
with \, $-\pi/2<\phi\leq \pi/2$ \,  it comes: 
\[
\theta_\perp
=\theta_1+\frac{\pi}{2}-\arctan \Big( 
\frac{(\alpha -\beta)\, \sin \theta_1\,  \cos \theta_1}
{\alpha\,\cos^2 \theta_1 + \beta \,\sin^2 \theta_1} \Big) 
\]
as well as: \quad $\dsm
\cos\theta_\perp=\frac{-\beta\, \sin \theta_1 }
{\sqrt{\alpha^2\, \cos^2 \theta_1 +\beta^2\, \sin^2 \theta_1}}
\quad , \quad 
\sin\theta_\perp=\frac{\alpha\, \cos \theta_1 }
{\sqrt{\alpha^2\, \cos^2 \theta_1+\beta^2\, \sin^2  \theta_1}}
$
\quad ,
\[
\cos \theta_1 =\frac{\beta\, \sin\theta_\perp }
{\sqrt{\alpha^2\, \cos^2\theta_\perp +\beta^2\, \sin^2\theta_\perp}}
\quad \hbox{and} \quad 
\sin \theta_1=\frac{-\alpha\, \cos\theta_\perp }
{\sqrt{\alpha^2\, \cos^2\theta_\perp +\beta^2\, \sin^2\theta_\perp}}
\quad .
\]

The main point here is that $\theta_{s+1}$ approximates indeed
$y(1-h)$ instead of $y(1)$\,, thus one more half-step of our middle-point method 
will give a better estimate of $y(1)$, hence a sharper boundary condition 
for the differential equation.
For each $\phi\in \MR$ we have:
\[
g(\phi)\,=\,
\frac{-\alpha\,\beta }
{\sqrt{\alpha^2\, \cos^2\theta_\perp +\beta^2\, \sin^2\theta_\perp}}
\,\sin\big(\phi-\theta_\perp)
\,=\,
-\frac{\sin\big(\phi-\theta_\perp)}
{f(\theta_\perp)}
\]
hence we obtain:
\, $
\theta_{s+1}
\,=\,
\theta_\perp-\arcsin \big( h\,f(\theta_\perp) \big)
$ \, , \, but we have: $\norme{f}_\infty\leq 3$ 
thus it comes:
\[
\absol{\theta_{s+1}-\big(\theta_\perp-h\, f(\theta_\perp)\big)}
\,\leq \,
\frac{1+2\,(3\,h)^2}{6\,\big(1-(3\,h)^2\big)^\frac{5}{2}}\,(3\,h)^3
\,\leq \,
4.7\, h^3
\]
since we assume that $s\geq 15$\,, hence $h\leq 1/31$\,. 
The function $z$ such that \, $z'=f(z)$ \, and \,
$z (1)=\theta_\perp$ \, also
satisfies: \, $\absol{z''(t)} \,\leq 4$ \, for all $t$ hence:
\[
\absol{z (1-h )-(\theta_\perp-h\, f(\theta_\perp))}
\, \leq \, 2\, h^2  \quad ,
\]
thus we obtain: \, 
$ \absol{z (1-h )-\theta_{s+1}}
\, \leq \, 2\, h^2 + 4.7\,h^3
$ \, , \,
hence by the above estimate:
\[
\absol{z (1-h )-y(1-h )}
\, \leq \, 2\, h^2   + 4.7\,h^3 +
5\,h^2 \, \big( e^{\,2\,(1-h)}-1\big)  
\quad .
\]
But $f$ is $2$-lipschitz, thus we get:
\, $
\absol{y(1)-\theta_\perp}
\,=\,
\absol{z(1)-y(1)}
\, \leq \,
 e^{\,2\,h}\,\absol{z(1-h)-y(1-h)}
$ \, and finally:
\, $
\absol{y(1)-\theta_\perp}
\, \leq \, 
e^{\,2\,h}\, \big(2\, h^2   + 4.7\,h^3 \big)
+ 5\,h^2 \, \big( e^{\,2}-e^{\,2\,h}\big)  
$ \, , \, hence:
\[
\Absole{y(1)-\theta_\perp}
\, \leq \, 
\big(5 \, e^{\,2}-(3- 4.7\,h)\,e^{\,2\,h}\big)\,h^2
\, \leq \, 33.95 \,h^2
 \, .
\]

Moreover, we obtain:
$
g(\phi)=
-\sqrt{\alpha^2\, \cos^2 \theta_1+\beta^2\, \sin  \theta_1}
\,\sin\big(\phi-\theta_\perp)$  for all $\phi\in \MR$\,,
hence:
\[
\theta_{s+1}
\,=\,
\theta_\perp
-\arcsin \Big(\frac{h }
{\sqrt{\,\alpha^2\, \cos^2\theta_1 +\beta^2\, \sin\theta_1\,}} 
\Big) \quad ,
\]
thus: \quad $\dsm
\Absole{\theta_{s+1} - \Big( \theta_\perp
-\frac{h }{\sqrt{\,\alpha^2\, \cos^2\theta_1 +\beta^2\, \sin\theta_1\,}}
\Big)}
\,\leq \, 4.7\,h^3
$ \quad as above, then it comes:
\[
\Absole{\sin\theta_{s+1} - \Big( \sin\theta_\perp
-\frac{h\,\cos\theta_\perp }
{\sqrt{\,\alpha^2\, \cos^2\theta_1 +\beta^2\, \sin\theta_1\,}}
\Big)}
\,\leq \, 4.7\,h^3 +\frac{1}{2}\,(3\,h)^2
\,\leq \, 4.66\,h^2 
\]
since \, $\alpha \geq \beta \geq 1/3$ \, , \, hence:
\[
\Absole{\sin\theta_{s+1} - \Big( 
\frac{\alpha\, \cos \theta_1 }
{\sqrt{\alpha^2\, \cos^2 \theta_1+\beta^2\, \sin  \theta_1}}
+\frac{\beta\,h\, \sin \theta_1 }
{(\alpha^2\, \cos^2 \theta_1 +\beta^2\, \sin \theta_1)}
\Big)}
\,\leq \, 4.66\,h^2 
\quad ,
\]
and we get the following estimate which will be useful in the next
section: 
\[
\Absole{ \sin \theta_1 \,\sin\theta_{s+1} - \Big( 
\frac{\alpha\, \sin \theta_1 \,\cos \theta_1 }
{\sqrt{\alpha^2\, \cos^2 \theta_1+\beta^2\, \sin  \theta_1}}
+\frac{\beta\,h\, \sin^2 \theta_1 }
{(\alpha^2\, \cos^2 \theta_1 +\beta^2\, \sin \theta_1)}
\Big)}
\,\leq \, 4.66\,h^2 
\quad .
\]

\section{The integral equation}\label{en}

The norm equation writes:
\[
2\,h\,\sum_{k=1}^{s} \cos\big(\,\theta_k+\theta_{k+1}\big)
\,=\,
\frac{(2\,s+2)\,\beta-2\,s\,\alpha}{(2\,s+1)\,(\alpha+\beta)}
-2\,\beta\,\sin \theta_{1}\, \sin \theta_{s+1}  
 \quad ,
\]
and for all $1\leq k \leq s$ we obtained in section \ref{equadiff}:
\[
\Absole{\frac{\theta_k+\theta_{k+1}}{2}-y (t_{2\,k-1})}
\, \leq \,
5\,h^2 \, \big( e^{\,2\,t_{2\,k-1}}-1\big)+2\,h^2 
\]
whenever $s\geq 15$\,, 
thus:
\[
\Absole{\cos\big(\theta_k+\theta_{k+1}\big)
-\cos\big(2\,y(t_{2\,k-1})\big)}
\, \leq \,
10\,h^2 \, \big( e^{\,(4\,k-2)\,h}-1\big)+4\,h^2 
\quad .
\]
If \,  $\dsm
I_s =2\,h\,\sum_{k=1}^{s}\cos\big(2\,y(t_{2\,k-1})\big)$  
denotes the approximation of 
$\dsm 
\int_0^{1-h}\cos\big(\,2\,y(t)\big)\,dt
$
obtained by the middle point method, we get first:
\[
\Delta_s=
\Absole{
2\,h\,\sum_{k=1}^{s} \cos\big(\,\theta_k+\theta_{k+1}\big)
-I_s}
\, \leq \,
2\,h\,\sum_{k=1}^{s}
\Big(\,10\,h^2 \, \big(\, e^{\,(4\,k-2)\,h}-1\,\big)+4\,h^2 \,\Big)
\]
\[
\, = \,
4\,h^3\,\Big(\,5\, e^{\,2\,h}\,\frac{ e^{\,4\,s\,h}-1}{e^{4\,h}-1}
-3\,s\,\Big) 
\, \leq \,
5\,h^2\, e^{\,2\,h}\,\big(e^{\,4\,s\,h}-1\big)-12\,s\,h^3
\, \leq \,
\big(5\,e^{\,2}-10.8\big)\,h^2
\quad .
\]
If we set now: \, $g(t)=\cos \big(\,2\,y(t)\,\big)$  \,
for all $t\in \MR$\,, we get:
\[
g''(t)=-4\, \cos \big(\,2\,y(t)\,\big)\,\big(\,y'(t)\,\big)^2
-2\, \sin \big(\,2\,y(t)\,\big)\,\,y''(t)
\]
\[
=  \Big(\frac{1}{\beta^2}-\frac{1}{\alpha^2}\Big)
-2\,\Big(\, \frac{1}{\alpha^2}+\frac{1}{\beta^2}\,\Big)\,
\cos \big(\,2\,y(t) \,\big)\,
-3\,\Big(\frac{1}{\beta^2}-\frac{1}{\alpha^2}\Big)\,
\cos^2 \big(\, 2\,y(t) \,\big)
\]
hence we obtain: \,
$\norme{g''}_\infty \,\leq \, 36$ \, , \, thus
the middle point method gives the estimate:
\[
\Absole{
\int_0^{1-h}\cos\big(\,2\,y(t)\big)\,dt
-I_s}
\, \leq \,
36\,\frac{\big(2\,h\big)^2}{24} \quad ,
\]
and we obtain: 
\[
\Absole{
2\,h\,\sum_{k=1}^{s} \cos\big(\,\theta_k+\theta_{k+1}\big)
-\int_0^{1-h}\cos\big(\,2\,y(t)\big)\,dt
}
\, \leq \,
\big(\,5\,e^{\,2}-4.8\,\big)\,h^2
\quad .
\]

But we also have: \,$\norme{g'}_\infty \,\leq \, 6$ \, , \, thus:
\[
\Absole{\int_0^{1}\cos\big(\,2\,y(t)\big)\,dt
-\int_0^{1-h}\cos\big(\,2\,y(t)\big)\,dt-
h\,\cos\big(\,2\,y(1)\big)}
\,\leq \,
3\,h^2
\]
hence: \, $\absol{I-S}\,\leq\, \big(\,5\,e^{\,2}-1.8\,\big)\,h^2$ \, where:
\[
I\,=\,\int_0^{1}\cos\big(\,2\,y(t)\big)\,dt
\quad \hbox{and} \quad
S\,=\,
2\,h\,\sum_{k=1}^{s} \cos\big(\,\theta_k+\theta_{k+1}\big)
+h\,\cos\big(\,2\,y(1)\,\big) \quad ,
\]
thus we get: \,
$\dsm
S\,=\,
\frac{(2\,s+2)\,\beta-2\,s\,\alpha}{(2\,s+1)\,(\alpha+\beta)}
-2\,\beta\,\sin \theta_{1}\, \sin \theta_{s+1}  
+h\,\cos\big(\,2\,y(1)\,\big)$ \,
by the norm equation. But we obtained above: 
\[
\Absole{ \sin \theta_1 \,\sin\theta_{s+1} - \Big( 
\frac{\alpha\, \sin \theta_1 \,\cos \theta_1 }
{\sqrt{\alpha^2\, \cos^2 \theta_1+\beta^2\, \sin  \theta_1}}
+\frac{\beta\,h\, \sin^2 \theta_1 }
{(\alpha^2\, \cos^2 \theta_1 +\beta^2\, \sin \theta_1)}
\Big)}
\,\leq \, 4.66\,h^2 
\]
and \,
$\dsm \absol{y(1)-\theta_\perp}
\,\, \leq \, 
\big(5 \, e^{\,2}-3 \big)\,h^2
$ \, ,
thus:
$\dsm \absol{h\,\cos\big(2\,y(1)\big)-h\,\cos\big(2\,\theta_\perp\big)}
\,\, \leq \, 
\big(10 \, e^{\,2}-6 \big)\,h^3
$\,,
where: \,
 $\dsm
\cos\theta_\perp=\frac{-\beta\, \sin \theta_1 }
{\sqrt{\alpha^2\, \cos^2 \theta_1 +\beta^2\, \sin^2 \theta_1}}
$ \, and \, $
\dsm \sin\theta_\perp=\frac{\alpha\, \cos \theta_1 }
{\sqrt{\alpha^2\, \cos^2 \theta_1+\beta^2\, \sin^2  \theta_1}}
$ \, , \, hence it comes:
$\dsm \absol{S-T}
\,\, \leq \, 
 4.66\,h^2 + \big(10 \, e^{\,2}-6 \big)\,h^3
$ \, where:
\[
T\,=\,
\frac{(2\,s+2)\,\beta-2\,s\,\alpha}{(2\,s+1)\,(\alpha+\beta)} 
-\frac{2\,\alpha\,\beta\,\sin \theta_{1}\,  \cos \theta_1 }
{\sqrt{\alpha^2\, \cos^2 \theta_1+\beta^2\, \sin^2  \theta_1}}
-\frac{2\,\beta^2\, \sin^2 \theta_1 }
{(2\,s+1)\,(\alpha^2\, \cos^2 \theta_1 +\beta^2\, \sin^2 \theta_1)}
\]
\[
+\,\,\,\frac{1}{2\,s+1}\,\,
\frac{\beta^2\, \sin^2 \theta_1 -\alpha^2\, \cos^2 \theta_1}
{\alpha^2\, \cos^2 \theta_1 +\beta^2\, \sin^2 \theta_1}
\]
\[
\,=\,
\frac{(2\,s+2)\,\beta-2\,s\,\alpha}{(2\,s+1)\,(\alpha+\beta)} 
-\frac{1}
{2\,s+1}
-\frac{2\,\alpha\,\beta\,\sin \theta_{1}\,  \cos \theta_1 }
{\sqrt{\alpha^2\, \cos^2 \theta_1+\beta^2\, \sin^2  \theta_1}}
\]
\[
\,=\,
\frac{\beta-\alpha}{\alpha+\beta} 
-\frac{2\,\alpha\,\beta\,\sin \theta_{1}\,  \cos \theta_1 }
{\sqrt{\alpha^2\, \cos^2 \theta_1+\beta^2\, \sin^2  \theta_1}}
\quad .
\]
We get finally the integral equation for all integers $s\geq 15$\,:
\[
\Absole{\int_0^{1}\cos\big(\,2\,y(t)\big)\,dt
+\frac{\alpha-\beta}{\alpha+\beta} 
+\frac{2\,\alpha\,\beta\,\sin \theta_{1}\,  \cos \theta_1 }
{\sqrt{\alpha^2\, \cos^2 \theta_1+\beta^2\, \sin^2  \theta_1}}
}
\,\, \leq \, 
C_1\,h^2
\]
where: \,
$\dsm C_1
\,=\,  4.66+ \big(10 \, e^{\,2}-6 \big)\,\frac{1}{31} +5\,e^{\,2}-1 .8
\, \leq \, 41.99$ \, , \,
but $f$ never vanishes, thus:
\[
I\,=\,
\int_0^{1}\cos\big(\,2\,y(t)\big)\,dt
\,=\,
\int_0^{1}\cos\big(\,2\,y(t)\big)\frac{y'(t)}{f\big(y(t)\big)}\,dt
\,=\,
\int_{\theta_1}^{y(1)}\frac{\cos \big(2\,x\big)}{f(x)}\,dx
\]
and the integral equation becomes a second boundary condition:
\[
\Absole{
\beta\,\int_{\theta_1}^{y(1)}\frac{\cos \big(2\,x\big)}
{\sqrt{\,1-\big(1-\frac{\beta^2}{\alpha^2}\big)\,\sin^2 x}}\,dx
+\frac{\alpha-\beta}{\alpha+\beta} 
+\frac{2\,\alpha\,\beta\,\sin \theta_{1}\,  \cos \theta_1 }
{\sqrt{\alpha^2\, \cos^2 \theta_1+\beta^2\, \sin  \theta_1}}
}
\,\leq\,
C_1\,h^2
\quad .
\]

\section{Numerical study of the asymptotic case} \label{numerique-ed} \label{bilan}

If $s\geq 15$ and the sequence $(\theta_k)_{1\leq k\leq s+1}$ 
from section \ref{equadiff} exists (that is, if 
$\lambda_2^{2s+1}>\lambda_2^{2s-1}$), we set 
$\theta=\theta_1\,\in[\, 0\vge \pi\,]$ and we proved that
the solution $y\dpe [\,0\vge 1]\dans \MR$ of: 
\[
y'=f(y)=\frac{1}{\beta}\,\sqrt{\,1
-\Big(1-\frac{\beta^2}{\alpha^2}\Big)\,\sin^2 y\,}
\]
such that \, $y(0)=\theta$ \, satisfies:
\[ 
\absol{CB(\alpha\vge \beta \vge \theta)}\,  
\,\leq\, \frac{33.95}{\big(2\,s+1\big)^2}
\quad \hbox{and}\quad
\absol{EN(\alpha\vge \beta \vge \theta)}\, 
\, \leq \,
\frac{41.99}{\big(2\,s+1\big)^2}
\quad ,
\]
where: 
\[
EN(\alpha\vge \beta \vge \theta)
\,=\,
\beta \,
\int_{\theta}^{\,y(1)}\frac{\cos \big(2\,x\big)}
{\sqrt{\,1-\big(1-\frac{\beta^2}{\alpha^2}\big)\,\sin^2 x}}\,dx
+\frac{\alpha-\beta}{\alpha+\beta} 
+\frac{2\,\alpha\,\beta\,\sin \theta\,  \cos \theta }
{\sqrt{\alpha^2\, \cos^2 \theta+\beta^2\, \sin  \theta}}
\]
and:
\quad 
$
CB(\alpha\vge \beta \vge \theta) \, = \, y(1)-\theta_\perp
$
\quad 
where
\quad 
$\dsm
\theta_\perp \, = \,
\theta+\frac{\pi}{2}-\arctan \Big( \,
\frac{(\alpha -\beta)\, \sin \theta\,  \cos \theta}
{ \alpha \, \cos^2 \theta + \beta \,\sin^2 \theta} \, \Big) 
$\,.
Setting: \,
$\dsm \gamma\,=\,\frac{\beta}{\alpha}\,\in \Big[\,\frac{1}{3}\vge 1\,\Big]$
\, the differential equation solves as: 
\[
\int_\theta^{y(t)}
\frac{dx}{\sqrt{\,1-(1-\gamma^2)\,\sin^2 x\,}\,} 
\, = \, 
\frac{t}{\beta}
\]
 for all $t\in \MR$\,, and we are left  with $3$ equations in $3$ variables 
(including a Jacobi function) which turn out to be incompatible.
In the limit case we get: \,
$CB(\alpha \vge \beta \vge \theta)\,=\,0$ \, , thus:
\[
y(1)
\,=\,
\theta+\frac{\pi}{2}-\arctan \Big(\, 
\frac{(1 -\gamma)\, \sin \theta\,  \cos \theta}
{\cos^2 \theta + \gamma \,\sin^2 \theta} \,\Big) 
\,=:\,\varphi(\gamma\vge \theta) \quad ,
\]
hence: 
$1/\beta\,=B\,(\gamma\vge \theta)$  and  $\alpha\,=\,\beta/\gamma$ 
where: \,
 $ \dsm
B(\gamma\vge \theta)
\,=\, \int_\theta^{\varphi(\gamma\vge \theta)}
\frac{dt}
{\,\sqrt{\,\cos^2 t +\gamma^2\,\sin^2 t\,}\,} 
$ \, ,
thus:
\, $\dsm
\frac{\partial B}{\partial \theta}(\gamma\vge \theta)
= 
\frac{1}{\sqrt{\,\cos^2 \theta_\perp+\gamma^2\,\sin^2 \theta_\perp}}\, 
\frac{\partial \varphi}{\partial \theta}(\gamma\vge \theta)
-\frac{1}{\sqrt{\,\cos^2 \theta +\gamma^2\,\sin^2 \theta\,}\,} 
$ \, where:
\[
\frac{\partial \varphi}{\partial \theta}(\gamma\vge \theta)
=\frac{\gamma }
{\cos^2 \theta + \gamma^2 \,\sin^2 \theta} 
\quad \hbox{and we obtained:} \quad
\sin\theta_\perp=\frac{\cos \theta }
{\sqrt{ \cos^2 \theta+\gamma^2\, \sin  \theta\,}\,}
\quad ,
\]
so we get finally:
\, $\dsm
\frac{\partial B}{\partial \theta}\,(\gamma\vge \theta)
\,=\,0
$ \, .
Therefore, the function $B$ doesn't depend on $\theta$\,, 
and if we choose $\theta=0$ we get: $\theta_\perp=\pi/2$\,, 
and we obtain a Jacobi function:
\[
B(\gamma\vge \theta)
\,=\, B(\gamma)
\,=\, \int_0^{\frac{\pi}{2}}
\frac{dt}
{\sqrt{\,\cos^2 t +\gamma^2\,\sin^2 t\,}\,} 
\,=\, 
InverseJacobiAM\,
\Big(\,\frac{\pi}{2}\vge \sqrt{\,1-\gamma^2\,}\, \Big)
\quad .
\]

We will
now compute two real numbers $\delta\vge \mu >0$ such that: 
\[
\absol{CB(\alpha\vge \beta    \vge \theta)}
\,\leq\, \delta
\quad \hbox{and} \quad
\alpha+\beta \geq \frac{4}{3}
\quad \Longrightarrow \quad
EN(\alpha\vge \beta    \vge \theta)
\,\geq\,
\mu \quad ,
\]
which will allow us to conclude for $s$ large enough.
Setting:
$\nu=CB(\alpha \vge \beta \vge \theta)\in [-\delta \vge \delta\,]$
\, we get: \, $y(1)=\varphi(\gamma \vge \theta) + \nu$ \,  
and \,
$\alpha(\gamma \vge \theta \vge \nu)
=\beta(\gamma \vge \theta \vge \nu)\,/\,\gamma$ \, 
where:
\[
\frac{1}{\beta(\gamma \vge \theta \vge \nu)}
\,=\,
\quad \int_\theta^{\varphi(\gamma \vge \theta)+\nu}
\frac{dt}
{\sqrt{\,\cos^2 t +\gamma^2\,\sin^2 t\,}\,} 
\quad ,
\]
let 
$ \dsm
E(\gamma \vge \theta\vge \nu)
=\frac{1-\gamma}{1+\gamma} 
+\frac{\beta(\gamma \vge \theta\vge \nu)\,\sin (2\,\theta) }
{\sqrt{\,\cos^2 \theta +\gamma^2\,\sin^2 \theta\,}\,} 
+\beta(\gamma \vge \theta\vge \nu)\,
\int_{\theta}^{\varphi(\gamma \vge \theta)+\nu}
\frac{\cos(2\, t)\,dt}
{\sqrt{\,\cos^2 t +\gamma^2\,\sin^2 t\,}\,} 
$
and let 
$\dsm \Delta \subset  \Big[ \,\frac{1}{3}\vge 1 \,\Big]
\times[ \,0\vge \pi \,]\times [ \,-\delta \vge \delta \,]$
be the compact set defined by:
\[
\beta (\gamma \vge \theta \vge \nu) 
 \,\leq \, \alpha (\gamma \vge \theta \vge \nu) 
 \,\leq \, 1
\quad \hbox{and} \quad
\alpha (\gamma \vge \theta \vge \nu)
+\beta (\gamma \vge \theta \vge \nu) 
 \,\geq \, \frac{4}{3} \quad .
\]
We will estimate 
$\dsm \mu (\delta)
\,=\, 
\min_{(\gamma \vge \theta \vge \nu) \,\in\, \Delta}\, 
E(\gamma \vge \theta\vge \nu)
$ 
using only the following obvious lemma.
\begin{lemma} Let $p \in \MN^\ast$, let $K \subset \MR^p$ be a 
convex compact
set, let $F\dpe K \dans \MR$ be a $C^1$ function and for all
\, $1\leq k \leq p$ \, let
\, $\dsm
M_k
\, = \,
\max_{x\,\in\, K}\, 
\Absole{\, \frac{\partial\, F}{\partial\, x_k} \, (x)\,}
$\,.
Let $\cale \subset K$ be a finite set 
and let $\delta \in \MR_+^p$ such that
for all $x\in K$\,, there exists  $y \in \cale$ satisfying:
$\absol{x_k-y_k}\, \leq \delta_k$ for each $1\leq k \leq p$\,.
Then we have: \quad
$\rule{0cm}{0.55cm}\dsm \min_{x\in K} F(x)
\,\,\, \geq \,\,\,
 \min_{y\in \cale} F(y)
\,\,\, - 
\sum_{1 \,\leq\, k \,\leq\, p } M_k\,\delta_k
$
\, .
\end{lemma}
Hereafter, the set
$\cale$ will be called a net with step $2\,\delta_k$ in
the $k$-th variable and the number 
$\dsm \sum_{k=1 }^p M_k\,\delta_k$ \, will be called
the uncertainty of this approximation of the minimum of $F$\,.
\smallskip

Notice that if the set 
$K$ is not convex, but 
the function $F$ is $C^1$ on its convex hull $\widehat{K}$\,, it
suffices to obtain the estimates $M_k$ on 
$\widehat{K}$\,, and that the same is true if we get these
estimates on a set $D$ such that 
for each $x\in K$\,, there exists  $y \in \cale$ and a sequence 
$(z_k)_{\, 0\,\leq\, k \,\leq\, p }$ satisfying: $z_0=x$\,,
$z_p=y$\,, $z_{k}-z_{k-1}$ proportional to the $k$-th vector 
of the canonical basis of $\MR^p$ (up to permutation) 
and $[z_{k-1} \vge z_{k}]\subset D$ for each
$1\,\leq\, k \,\leq\, p$\,.

\bigskip

In order to conclude for $s\geq 15$ it suffices to prove that: \,
$\mu (0.0354) \geq 0.0446$ \, , \, 
but if we estimate roughly the derivatives of $E$ the uncertainty will be too 
big to conclude in a reasonable computation time: the worse case is when 
$\gamma$ is small because of all the terms in 
$1/\gamma^k$\,, and  we will first reduce the range in $\gamma$\,. We have:
\[
\frac{1}{\beta(\gamma \vge \theta \vge \nu)}
\,=\,
\int_{\theta}^{\varphi(\gamma,\theta)}
\frac{dt}{\sqrt{\cos^2 t+\gamma^2\,\sin^2 t}} 
+\int_{\varphi(\gamma,\theta)}^{\varphi(\gamma,\theta)+\nu}
\frac{dt}{\sqrt{\cos^2 t+\gamma^2\,\sin^2 t}} 
\]
hence : \, $ 
1/\beta(\gamma \vge \theta \vge \nu) 
\leq B(\gamma)+ \delta/\gamma  $ \, , \,
thus the condition: \, $\alpha=\beta/\gamma \leq 1$ \, 
implies:
\, $ \gamma\,B(\,\gamma) -1 \,\geq \,-\delta$ \,  where \,
$\delta\,=\,0.0354$ \,. \,
We get easily:
\[
 -\,\frac{\pi}{2\,\gamma^2}  \,\leq \, B'(\gamma)\,\leq\, 0
  \,\leq \, B'(\gamma)\,\leq\,\frac{\pi}{2\,\gamma} 
\quad \hbox{hence} \quad 
\Absole{\gamma\,B'(\gamma)+B(\gamma)}\,\leq\, 
\frac{\pi}{2\,\gamma}
\,<\, 5
\]
for all $\gamma \in [\,\,1/\,3\,\vge 1\,]$\,,
thus the Maple procedure:
\begin{verbatim}
n:=500; MM:=-1000:
for igammaa from 0 to n do gammaa:=evalf(1/3+(0.414-1/3)*igammaa/n): 
h:=gammaa*InverseJacobiAM(Pi/2,sqrt(1-gammaa^2))-1:
if h>MM then MM:=h: fi: od:  M:=evalf(MM + 5*(0.414-1/3)/(2*n));
\end{verbatim} 
answers: \, 
$\max \big\{\gamma\,B(\gamma)-1\vge 1/\,3\, \leq \gamma \leq 0.414 \big\}
\,<\,  -0.0355 \,<\, -\delta$ \, , \,  
hence we can assume that: 
\,  $0.414 \leq \gamma \leq 1$ \, . \, 
In order to estimate the
minimum of the function $E$\,, we have to estimate its partial
derivatives, which will be done in section \ref{deriveesEN}
by lengthy computations and numerical studies: we will
obtain the estimates:
\[
\Absole{\frac{\partial E}{\partial \nu}(\gamma \vge \theta \vge \nu) }
\, \leq \,
2.48
 \quad , \quad 
\Absole{\frac{\partial E}{\partial \theta}(\gamma \vge \theta \vge \nu) }
\, \leq \,
4.41
\quad \hbox{and} \quad
\Absole{\frac{\partial E}{\partial \gamma}(\gamma \vge \theta \vge \nu) }
\, \leq \,
4.33
\]
for all 
$(\gamma \vge \theta \vge \nu)\in D= [\, 0.414 \vge 1 \,]
\times[ \,0\vge \pi \,]\times [ \,-0.0354 \vge 0.0354 \,]$\,.
The Maple procedure from appendix \ref{procmu} 
realizes an uncertainty of \, $0.1514/n$ \, 
on the minimum $\mu (\delta)$ of the function
$E$ on the domain $\Delta \subset D$ where $D$ is convex.
If $n=6$, it answers in  60 minutes:  
\, $\mu(0.0354) \,>\, 0.0484 \,>\, 0.0446$ \, , \,
and we can conclude that: \, 
$\lambda_{\,2}^{2\,s+1}\, = \,\lambda_{\,2}^{2\,s-1}$ \,
for all integers \, $s\geq 15$ \, .

\section{Iteration in the initial cases}\label{iteration}

To complete the proof, it remains to deal with the cases where $2\leq s\leq 14$\,, 
and to keep the computation time reasonable we will have to implement the
equations in C after simplifying them. 
First of all, we rewrite the induction relation as obtained in 
section \ref{essentiels}:
\[
\alpha^2\,\cos^2\Big(\,\frac{\theta_{k}+\theta_{k+1}}{2}\,\Big) 
+\beta^2\,\sin^2\Big(\,\frac{\theta_{k}+\theta_{k+1}}{2}\,\Big) 
\,=\,
(2s+1)^2\,\alpha^2\,\beta^2 \,
\sin^2\Big(\,\frac{\theta_{k}-\theta_{k+1}}{2}\,\Big) 
\]
to get: \, $
\alpha^2\, \big(1+\cos (\theta_{k}+\theta_{k+1})\big)+
\beta^2 \, \big(1-\cos (\theta_{k}+\theta_{k+1})\big)
\,=\,
(2s+1)^2\,\alpha^2\,\beta^2 \,
\big( 1-\cos (\theta_{k}-\theta_{k+1}) \big)$ 
\,  hence: \, $
\big( (2s+1)^2\,\alpha^2\,\beta^2 -\alpha^2+\beta^2\big) \,
\sin\theta_{k}\,\sin \theta_{k+1}
+
\big( (2s+1)^2\,\alpha^2\,\beta^2 +\alpha^2-\beta^2\big) \,
\cos\theta_{k}\,\cos \theta_{k+1}
$ \, 
\[
= \quad 
(2s+1)^2\,\alpha^2\,\beta^2 -\alpha^2-\beta^2
\quad , 
\]
that is: \quad $
\big(1-2\,h^2\,v\big)\,\sin\theta_{k}\,\sin \theta_{k+1}
+
\big(1+2\,h^2\,v\big)\,\cos\theta_{k}\,\cos \theta_{k+1}
=
1-2\,h^2\,u
$ \quad where:
\[
h=\frac{1}{2\,s+1}
\quad , \quad 
u=\frac{1}{2}\,\Big(\frac{1}{\beta^2}+\frac{1}{\alpha^2}\Big)
\quad \hbox{and} \quad 
v=\frac{1}{2}\,\Big(\frac{1}{\beta^2}-\frac{1}{\alpha^2} \Big)
\quad .
\]
This change of variables writes:
\, $ \dsm
\alpha (u \vge v)=\frac{1}{\sqrt{u-v}}
$ 
\, and
\, $ \dsm
\beta  (u \vge v)=\frac{1}{\sqrt{u+v}}
$ \, , \, 
where $(u\vge v)$  fulfill the inequalities: \,
$\dsm
1 \,\leq \, u \,\leq \,
5
$
\, , \,
$\dsm
0 \,\leq \, v \,\leq \, u-1 \,\leq \, 4
$
\,
and if \, $u\geq \dsm \frac{9}{4}$ \, :
\[
v 
\,\geq \, 
\sqrt{\frac{128\,u^2-144\,u-81-27\sqrt{9+32\,u}}{128}\,\,}
\quad ,
\]
and $D\subset \MR^2$ will hereafter denote the domain so defined.
If we set moreover:
\[
A=1+2\,h^2\,v
\quad ,\quad 
B=1-2\,h^2\,v
\quad \hbox{and} \quad 
C=1-2\,h^2\,u
\quad,
\]
then for all $1\leq k \leq s$ we have:  \quad $A\,\cos \theta_k\,\cos\theta_{k+1}
+
B\,\sin \theta_k\,\sin\theta_{k+1}
=C$ \quad  with the boundary condition: \,
$\dsm
\alpha\,   \cos \theta_1  \,  \cos \theta_{s+1} 
+ \beta\,   \sin \theta_{1}\, \sin \theta_{s+1} 
=\frac{1}{2s+1}
$ \,
and the norm equation:
\[
(2s+1)\,\beta\,\sin \theta_{1}\, \sin \theta_{s+1}  +
\sum_{k=1}^{s} \big(\cos\theta_k\,\cos\theta_{k+1} 
-\sin\theta_k\,\sin\theta_{k+1}\big)
\,=\,
\frac{(s+1)\,\beta-s\,\alpha}{\alpha+\beta} \quad .
\]

Setting:
\, $x_k=A\,\cos\theta_{k}$ \, and \, $y_k=B\,\sin\theta_{k}$ \,
we get:
\, $B^2\,x_k^2+A^2\,y_k^2=A^2\,B^2$ \, and:
\[
B\,x_k\,x_{k+1}+A\,y_k\,y_{k+1}=A\,B\,C
\]
for all $1\leq k \leq s$\,,
and after a few computations we get for all  $1\leq k\leq s$\,:
\[
x_{k+1} = A\,\frac{C\,x_{k}-y_{k}
\sqrt{x^2_{k}+y^2_{k}-C^2}}
{x^2_{k}+y^2_{k}} 
\quad \hbox{and} \quad 
y_{k+1} = B\,\frac{C\,y_{k}+x_{k}
\sqrt{x^2_{k}+y^2_{k}-C^2}}
{x^2_{k}+y^2_{k}} 
\quad .
\]
The boundary condition becomes:
\[
cb(u\vge v\vge \theta_1)=
\frac{\alpha}{A}\,   \cos \theta_1  \,  x_{s+1}  + 
\frac{\beta}{B}\,    \sin \theta_1  \, y_{s+1}
-\frac{1}{2\,s+1}=0
\]
and the norm equation becomes:
\[
en(u\vge v\vge \theta_1)=
\frac{s\,\alpha-(s+1)\,\beta}{(\alpha+\beta)\,(2s+1)}  
+
\frac{1}{2\,s+1}
\sum_{k=1}^{s} \Big(\frac {x_k\,x_{k+1}}{A^2}-
\frac {y_k\,y_{k+1}}{B^2}\Big)
+\beta\,\frac {\sin \theta_1 \, y_{s+1}}{B}
=0 \quad ,
\]
and we will now prove numerically that the function 
$m=\min (\absol{cb} \vge \absol{en})$ is nonzero, which will complete the proof.
The maxima of the derivatives of these iterate functions of 3 variables will be estimated
in section \ref{der-iteration}: the uncertainty on $m$ is at most the maximum of uncertainties on 
$cb$ and $en$ and the relevant domain $K=D\times [0\vge \pi]$ is not convex because the function
\[
u \mapsto
\sqrt{\frac{128\,u^2-144\,u-81-27\sqrt{9+32\,u}}{128}\,\,}
\]
is concave. However, going from $(u_0\vge v_0)$ to
$(u_1\vge v_1)$ by the segment
\, $(u_0\vge v_0)\rightarrow (u_0\vge v_1)$ \,
then by the segment: 
\, $(u_0\vge v_1)\rightarrow (u_1\vge v_1)$ \, we stay in $D$\,, hence it suffices to 
estimate these derivatives on $K$\,, and on this aim to estimate their second derivatives on $K$\,.
If $4\leq s \leq 14$ we will obtain in section \ref{der-iteration}:
\[
\Absole{\frac{\partial cb}{\partial \theta}}
\,\leq\,
16.5
\quad
,
\quad
\Absole{\frac{\partial cb}{\partial u}}
\,\leq\,
3.5
\quad
\hbox{and}
\quad
\Absole{\frac{\partial cb}{\partial v}}
\,\leq\,
3.6
\]
\[
\Absole{\frac{\partial en}{\partial \theta}}
\,\leq\,
24.1
\quad
,
\quad
\Absole{\frac{\partial en}{\partial u}}
\,\leq\,
5.15
\quad
\hbox{and}
\quad
\Absole{\frac{\partial en}{\partial v}}
\,\leq\,
6.25
\quad ,
\]
if $s=3$\,, we will get:
\[
\Absole{\frac{\partial cb}{\partial \theta}}
\,\leq\,
22.2
\quad
,
\quad
\Absole{\frac{\partial cb}{\partial u}}
\,\leq\,
4.35
\quad
\hbox{and}
\quad
\Absole{\frac{\partial cb}{\partial v}}
\,\leq\,
4.4
\]
\[
\Absole{\frac{\partial en}{\partial \theta}}
\,\leq\,
32.9
\quad
,
\quad
\Absole{\frac{\partial en}{\partial u}}
\,\leq\,
6.46
\quad
\hbox{and}
\quad
\Absole{\frac{\partial en}{\partial v}}
\,\leq\,
6.9
\quad ,
\]
and the worst case will be $s=2$\,, where the minimum of $m$ is the smallest 
hence we need sharper estimates, and we will obtain steadily:
\[
\Absole{\frac{\partial cb}{\partial \theta}}
\,\leq\,
7.56
\quad
,
\quad
\Absole{\frac{\partial cb}{\partial u}}
\,\leq\,
1.39
\quad
\hbox{and}
\quad
\Absole{\frac{\partial cb}{\partial v}}
\,\leq\,
1.39
\]
\[
\Absole{\frac{\partial en}{\partial \theta}}
\,\leq\,
10.09
\quad
,
\quad
\Absole{\frac{\partial en}{\partial u}}
\,\leq\,
2.22
\quad
\hbox{and}
\quad
\Absole{\frac{\partial en}{\partial v}}
\,\leq\,
2.81
\quad .
\]
The C procedure from appendix \ref{minimumm}  
estimates the minimum of
$m$ with a step of $1/(4\,n)$
as well as the uncertainty $\delta m$\,: 
for $4 \leq s \leq 14$ it requires $n=121$
to answer in 2 hours: $\min m= 0.0367... > 0.0365... = \delta m$ 
thus $\min m >0$,
for $s=3$ it needs $n=200$ and answers in 25 minutes: 
$\min m= 0.0339... > 0.0221... = \delta m$
and for $s=2$ it requires $n=200$ and answers
$\min m= 0.02725... > 0.02205... = \delta m$ in 25 minutes. The
Gr{\"u}nbaum conjecture is thus proved up to the estimates of the
derivatives in the final 2 sections. 


\section{Estimates on the partial derivatives of $E$}
\label{deriveesEN}

On the convex domain \,
$[\, 0.414 \vge 1 \,]
\times[ \,0\vge \pi \,]\times [ \,-\delta \vge \delta \,]$ \, 
where \,  $\delta=0.0354$ \, 
we have:
\[
E(\gamma \vge \theta\vge \nu)
=\frac{1-\gamma}{1+\gamma} 
+\beta(\gamma \vge \theta\vge \nu)\,
\Big(\,\frac{\sin (2\,\theta) }
{\sqrt{\,\cos^2 \theta +\gamma^2\,\sin^2 \theta}} 
+\int_{\theta}^{\,y(1)}
\frac{\cos(2\, t)\,dt}{\sqrt{\,\cos^2 t +\gamma^2\,\sin^2 t}} 
\,\Big)
\]
where:  \, $y(1)=\theta_\perp+\nu$ \, , \,
\, $\dsm
\theta_\perp
\,=\,
\theta+\frac{\pi}{2}-\arctan \Big(\, 
\frac{(1 -\gamma)\, \sin \theta\,  \cos \theta}
{\cos^2 \theta + \gamma \,\sin^2 \theta} \,\Big) 
$ \, and:
\[
\frac{1}{\beta(\gamma \vge \theta \vge \nu)}
\,=\,
\int_{\theta}^{\,y(1)}
\frac{dt}{\sqrt{1-\big(1-\gamma^2\big)\,\sin^2 t}} 
\,=\,
B(\gamma)
+\int_{\varphi(\gamma,\theta)}^{\varphi(\gamma,\theta)+\nu}
\frac{dt}{\sqrt{\cos^2 t+\gamma^2\,\sin^2 t}} 
\quad ,
\]
thus: \, 
$\theta \leq \theta_\perp \leq \theta+\pi$ \, . \,
In order to estimate the partial derivatives of $E$\,, 
the most efficient way
is a numerical study of their approximate expressions of 2 variables, where we get
rid of the small variable $\nu$\,.
First we get, since the function $B$ decreases:
\[
1.485
\leq
B(1)-\frac{\delta}{\gamma}
\leq
B(\gamma)-\frac{\delta}{\gamma}
\leq
\frac{1}{\beta(\gamma \vge \theta \vge \nu)}
\leq
B(\gamma)+\frac{\delta}{\gamma}
\leq
B(0.414)+\frac{\delta}{\gamma}
\leq 2.414
\]
hence:  $0.4143 \leq \beta_{min} \leq \beta \leq \beta_{max} \leq 0.6733$\,, and moreover:
\[
\Absole{\beta(\gamma \vge \theta \vge \nu)-\frac{1}{B(\gamma)}}
\, \leq \, 
\frac{\delta}{B(\gamma)\,\big(\gamma\,B(\gamma)-\delta\big)}
\,=\,
\frac{\delta}{D(\gamma)}
\]
and we obtain: \,
$\absol{D'(\gamma)}\,\leq\pi/(2\,\gamma^2)\leq 29$\,, 
thus the Maple procedure:
\begin{verbatim}
n:=50000: mm:=100:
for igamma from 0 to n do gammaa:=evalf(0.414+0.586*igamma/n): 
B:=InverseJacobiAM(Pi/2,sqrt(1-gammaa^2)):h:=B*(gammaa*B-0.0354):
if h<=mm then mm:=h: fi: od:  m:=evalf(mm-29*0.586/(2*n)); 
dbeta:=0.0354/m; betamax:=evalf(2/Pi+dbeta);
\end{verbatim} 
shows that we have: \, 
$D(\gamma)\, > \, 2.16
$ \, ,  \, hence 
\, for all $(\gamma \vge \theta \vge \nu)\in D$ \, :
\[
\Absole{\beta(\gamma \vge \theta \vge \nu)-\widetilde{\beta}(\gamma)}
\, \leq \, 
\Delta_\beta
\,=\,
0.0164
\quad \hbox{where:}\quad 
\widetilde{\beta}(\gamma)=\frac{1}{B(\gamma)}
\quad ,
\]
and since the function $B$ decreases, we infer:
\[
0.4143
\, \leq \, 
\beta(\gamma \vge \theta \vge \nu) 
\, \leq \, \widetilde{\beta} (1) + \Delta_\beta
\,=\, \beta_{max}
\,\leq\, 0.6531 \quad .
\]

In what follows, we will repeatedly have to estimate the maximum 
on $[\, 0.414 \vge 1 \,] \times [ \,0\vge \pi \,]$ 
of functions similar to:
\[
f_{0} \dpe (\gamma \vge \theta) \mapsto 
\Absole{\frac{\sin \big(2\,\theta\big)}
{\sqrt{\,\cos^2 \theta+\gamma^2\,\sin^2 \theta}}}
\, \leq \, \frac{1}{\gamma}
\, \leq \, 2.42
\quad ,
\]
but this is a too rough estimate for us, and the study of $f_0 (\gamma \vge \cdot)$ 
shows that this maximum is indeed:
\[
M_0\,=\, \max_{0.414 \,\leq\, \gamma \,\leq\, 1} 
\Big(\frac{2}{1+\gamma}\Big)
\, \leq \, 1.4145
\quad ,
\]
but to spare us this lengthy study we first remark that 
the maximum in $\theta$ decreases with $\gamma$  hence 
its maximum (in $\gamma$)
is attained for $\gamma=0.414$\,, that we easily obtain:
\[
\Absole{\frac{\partial f_0}{\partial \theta}\,(\gamma\vge \theta)} 
\,\leq\,
\frac{2}{\gamma}+\frac{1}{2\,\gamma^3}
\quad ,
\]
hence the Maple procedure which will be thereafter denoted by 
$\star_{Maple}$\,:
\begin{verbatim}
gammaa:=0.414: derf0:=2/gammaa+1/(2*gammaa^3): n:=500000: MM0:=-1000:
for itheta from 0 to n do theta:=evalf(Pi*itheta/(n)):
f0:=evalf(abs(sin(2*theta))/(sqrt(1-(1-gammaa^2)*(sin(theta))^2))): 
if f0>MM0 then MM0:=f0: fi: od:  M0:=evalf(MM0+derf0*Pi/(2*n)); 
\end{verbatim} 
shows that \, $M_0 \,\leq\, 1.4145 $.
In order to simplify the notations, we will hereafter write:
\[
F(\gamma \vge \theta\vge \nu) \,=\,
G(\gamma \vge \theta\vge \nu)\,+(\leq)\,
 H(\gamma \vge \theta\vge \nu)
\]
instead of: \,
$\absol{F(\gamma \vge \theta\vge \nu)-G(\gamma \vge \theta\vge \nu)}
\, \leq \,\, 
\absol{H(\gamma \vge \theta\vge \nu)}$ \,
for all \,  $(\gamma \vge \theta\vge \nu)\, \in\, D$\,,
including the case where $G=0$\,. 
\bigskip

First of all, we have: \,
$y(1)\,=\,\theta_{\perp}\,+(\leq) \,\,\delta $ \,  and: 
\[
\frac{\partial E}{\partial \nu} =
\frac{\partial \beta}{\partial \nu}\,
\Big( \frac{\sin (2\,\theta) }
{\sqrt{\cos^2  \theta +\gamma^2\, \sin^2  \theta}}
+
\int_{\theta}^{y(1)}\frac{\cos(2\, t)\,dt}
{\sqrt{\,\cos^2 t+\gamma^2\,\sin^2 t\,}} \Big)
+
\frac{\beta\,\cos(2\,y(1) )}
{\sqrt{\,\cos^2 y(1)+\gamma^2\,\sin^2 y(1)}}
\]
where:
\, $\dsm
\frac{\partial \beta}{\partial \nu}=
\frac{-\beta^2}
{\sqrt{\,1-(1-\gamma^2)\,\sin^2 y(1)}}
=
\frac{-\beta^2}
{\sqrt{\,\cos^2 \theta_\perp+\gamma^2\,\sin^2 \theta_\perp}}
\,+ \leq \Big(\,\beta^2\,(1-\gamma^2)\,M_{1}\,\delta\,\Big)
$ \, and $M_{1}$ is the maximum of
\, $\dsm f_{1} \dpe (\gamma \vge \theta) \mapsto 
\frac{\absol{\sin \big(2\,\theta \big)}}
{\sqrt{\,\cos^2 \theta+\gamma^2\,\sin^2 \theta}^{\,3}}
$ \, , \, 
hence the above procedure $\star_{ Maple}$ 
shows that $M_{1} \,\leq\, 2.3471 $\,, thus:
\[
\Absole{\frac{\partial \beta}{\partial \nu}-
\frac{-\beta^2}
{\sqrt{\,\cos^2 \theta_\perp+\gamma^2\,\sin^2 \theta_\perp}}}
 \,\leq\, 
\Delta \beta_\nu
\,=\,
0.02945 \quad .
\]  
We also have: \, 
$\dsm
\frac{\cos(2\,y(1) )}{\sqrt{\,\cos^2 y(1)+\gamma^2\,\sin^2 y(1)}}
=
\frac{\cos(2\,\theta_\perp )}
{\sqrt{\,\cos^2 \theta_\perp +\gamma^2\,\sin^2 \theta_\perp}}
+\,  \leq \big(\,8\,M_2\,\delta\,\big)
$ 
\, , \, 
where
$M_2$ is the maximum of
\, $\dsm f_{\,2} \dpe (\gamma \vge\theta) \mapsto 
\frac{\sin \theta \, \cos^3 \theta}
{\sqrt{\,\cos^2 \theta+\gamma^2\,\sin^2 \theta}^{\,3}}
$ \, , \, 
and we obtain:
$\norme{f_{\,2}'}_\infty \,\leq\, 250$\,,
hence the Maple procedure $\star_{Maple}$ shows that \, $M_2\,\leq\, 0.56$ \, and:
\[
\Absole{
\frac{\cos(2\,y(1) )}{\sqrt{\,\cos^2 y(1)+\gamma^2\,\sin^2 y(1)}}
-\frac{\cos(2\,\theta_\perp )}
{\sqrt{\,\cos^2 \theta_\perp +\gamma^2\,\sin^2 \theta_\perp}}}
\,\leq\,
\Delta_{y(1)}
\,=\,
0.1649
\quad .
\] 
We also have:
\, $\dsm
\Absole{\frac{\partial \beta}{\partial \nu}}
\,\leq\, 
\frac{\beta^2}{\gamma}
\,=\,
M_{\beta_\nu}
\,\leq\, 
1.03
$ \, , \,  thus if we set:
\[
\overline{E_\nu}
=
\frac{-\beta^2}
{\sqrt{\cos^2\theta_\perp +\gamma^2\,\sin^2 \theta_\perp}}\,
\Big(\frac{\sin (2\,\theta) }
{\sqrt{\cos^2  \theta +\gamma^2\, \sin^2  \theta}}
+\int_{\theta}^{\theta_\perp}\frac{\cos(2\, t)\,dt}
{\sqrt{\,\cos^2 t+\gamma^2\,\sin^2 t\,}}\Big)
\]
\[
+\quad \frac{\beta\,\cos(2\,\theta_\perp )}
{\sqrt{\,\cos^2 \theta_\perp +\gamma^2\,\sin^2 \theta_\perp}}
\]
we get: \,
$\dsm
\Absole{
\frac{\partial E}{\partial \nu} -\overline{E_\nu} }
\, \leq \,
M_{\beta_\nu}\,\frac{\delta}{\gamma}
+{\Delta_\beta}\,\,\frac{1+\pi}{\gamma}
+\beta\,\Delta_{\, y(1)}
\, \leq \, 0.4904
$ \, , \, then if we set: 
\[
\widetilde{E_\nu} (\gamma \vge \theta)
=
\frac{-\big(\,\widetilde{\beta}(\gamma)\,\big)^{2}}
{\sqrt{\cos^2\theta_\perp +\gamma^2\,\sin^2 \theta_\perp}}\,
\Big(\frac{\sin (2\,\theta) }
{\sqrt{\cos^2  \theta +\gamma^2\, \sin^2  \theta}}
+\int_{\theta}^{\theta_\perp}\frac{\cos(2\, t)\,dt}
{\sqrt{\,\cos^2 t+\gamma^2\,\sin^2 t\,}}\Big)
\]
\[
+\quad \frac{\widetilde{\beta}(\gamma)\,\cos(2\,\theta_\perp )}
{\sqrt{\,\cos^2 \theta_\perp +\gamma^2\,\sin^2 \theta_\perp}}
\]
it comes:
$\dsm
\Absole{\overline{E_\nu}-\widetilde{E_\nu}}\,\leq\,
\frac{\,2\,\beta_{max}\,{\Delta_\beta}\,}{\gamma}\,\,
\frac{\,1+\pi\,}{\gamma}
+
\frac{\,{\Delta_\beta}\,}{\gamma}
\, \leq \, 0.5572
$
\, , \,
and finally:
\[
\Absole{
\frac{\partial E}{\partial \nu} (\gamma \vge \theta\vge \nu)
-\widetilde{E_\nu} (\gamma \vge \theta)}
\, \leq \,
1.048
\,\, .
\]

Similarly, we get:
\, $ \dsm
\frac{\partial E}{\partial \theta}=
\frac{\partial \beta}{\partial \theta}\,\Big(\,
\frac{\sin (2\,\theta) }
{\sqrt{\cos^2  \theta +\gamma^2\, \sin^2  \theta}}
+\int_{\theta}^{y(1)}\frac{\cos(2\, t)\,dt}
{\sqrt{\,\cos^2 t+\gamma^2\,\sin^2 t\,}}
\,\Big)
$
\[
+ \,
\frac{\beta\,(1-\gamma^2)\,\sin^2(2\,\theta )} 
{2\,\sqrt{\cos^2  \theta +\gamma^2\, \sin^2  \theta}^{\,\,3}}
\, + \,
\frac{\beta\,\cos(2\,y(1) )}
{\sqrt{\,\cos^2 y(1)+\gamma^2\,\sin^2 y(1)}} \, \,
\frac{\partial \varphi}{\partial \theta} 
\, - \,
\frac{\beta\,\cos(2\,\theta )}
{\sqrt{\,\cos^2 \theta+\gamma^2\,\sin^2 \theta}} 
\quad ,
\]
where:
$\dsm 
\frac{1}{\beta(\gamma \vge \theta \vge \nu)}
\,=\,
B(\gamma)
+
\int_{\varphi(\gamma,\theta)}^{\,\varphi(\gamma,\theta)+\nu}
f_{\,3}(t)\,dt
$ \, and \,
$\dsm f_{\,3}(t)=
\frac{1}{\sqrt{\cos^2 t+\gamma^2\,\sin^2 t}} 
$\,, hence:
\[
\frac{-1}{\beta^2}\, \frac{\partial \beta}{\partial \theta}
\,=\,
\big( f_{\,3}(\theta_\perp+\nu)- f_{\,3}(\theta_\perp)\big)\,
\frac{\partial \varphi}{\partial \theta}
\,=\,
\nu\, f_{\,3}'(\theta_\perp+s\,\nu)\,
\frac{\partial \varphi}{\partial \theta} 
\,=(\leq)\, \Big(\,(1-\gamma^2)\,M_{\,1}\,\delta\,\Big)\,
\frac{\partial \varphi}{\partial \theta} 
\]
where $0<s<1$
since $f_{3}'=(1-\gamma^2)\,f_{\,1}$\,,
and: \,  $ \dsm
\frac{\partial \varphi}{\partial \theta}
=\frac{\gamma }
{\cos^2 \theta + \gamma^2 \,\sin^2 \theta}
\,=(\leq)\, \Big(\frac{1}{\gamma}\Big)
$ \, , \, 
thus if we set:
\[
\overline{E_\theta}
\, = \, 
\frac{\beta\,(1-\gamma^2)\,\sin^2(2\,\theta )} 
{2\,\sqrt{\cos^2  \theta +\gamma^2\, \sin^2  \theta}^{\,\,3}}
\, + \,
\frac{\beta\,\cos(2\,\theta_\perp )}
{\sqrt{\,\cos^2 \theta_\perp+\gamma^2\,\sin^2 \theta_\perp}} \, \,
\frac{\partial \varphi}{\partial \theta}
\, - \,
\frac{\beta\,\cos(2\,\theta )}
{\sqrt{\,\cos^2 \theta+\gamma^2\,\sin^2 \theta}} 
\]
it comes: 
$\dsm
\Absole{\frac{\partial E}{\partial \theta}-\overline{E_\theta}}
\,\leq \,
\frac{(1+\pi+\delta)\,\beta_{max}^2\,(1-\gamma^2)\,M_{\,1}\,\delta}
{\gamma^2}
\,+\,\frac{\beta_{max}\,\Delta_{y(1)}}{\gamma}
\, < \, 1.3437
$\,,  hence if: 
\[
\widetilde{E_\theta}
\, = \,
\frac{\widetilde{\beta}\,(1-\gamma^2)\,\sin^2(2\,\theta )} 
{2\,\sqrt{\cos^2  \theta +\gamma^2\, \sin^2  \theta}^{\,\,3}}
\, + \,
\frac{\widetilde{\beta}\,\cos(2\,\theta_\perp )}
{\sqrt{\,\cos^2 \theta_\perp+\gamma^2\,\sin^2 \theta_\perp}} \, \,
\frac{\partial \varphi}{\partial \theta}
\, - \,
\frac{\widetilde{\beta}\,\cos(2\,\theta )}
{\sqrt{\,\cos^2 \theta+\gamma^2\,\sin^2 \theta}} 
\]
we get: \,
$\dsm
\Absole{\overline{E_\theta}-\widetilde{E_\theta}}
\,\leq \,
\Big(\frac{(1-\gamma^2)}{2\,\gamma^3}
\,+\,\frac{1}{\gamma^2}
\,+\,\frac{1}{\gamma} \,\Big)\,
{\Delta_\beta}
\, < \, 0.1919
$
\, , \, and finally:
\[
\Absole{\frac{\partial E}{\partial \theta}(\gamma \vge \theta\vge \nu)-\widetilde{E_\theta}(\gamma \vge \theta)}
\,\leq \,
1.5356 \quad .
\]

\noindent
We also have: \, $\dsm 
\frac{\partial E}{\partial \gamma}
\,=\,
\frac{\partial \beta}{\partial \gamma}\, 
\Big(\, \frac{\sin (2\,\theta) }
{\sqrt{\cos^2  \theta +\gamma^2\, \sin^2  \theta}}
\,+\,\int_{\theta}^{y(1)}\frac{\cos(2\, t)\,dt}
{\sqrt{\,\cos^2 t+\gamma^2\,\sin^2 t\,}}
\,\Big)
-\frac{2}{\big( 1+\gamma \big)^2} 
$
\[
-\gamma\,\beta\,\Big(
\frac{\sin (2\,\theta)\, \sin^2  \theta }
{\sqrt{\cos^2  \theta +\gamma^2\, \sin^2  \theta}^{\,\,3}}
\,+\,
\int_{\theta}^{y(1)}\frac{\cos(2\, t)\,\sin^2 t\,dt}
{\sqrt{\,\cos^2 t+\gamma^2\,\sin^2 t\,}^{\,\,3}}
\Big)
+
\frac{\beta\,\cos\big(2\, y(1)\big)}
{\sqrt{\cos^2 y(1)+\gamma^2\,\sin^2 y(1) }}\,
\frac{\partial \varphi}{\partial \gamma}
\]
where:
\, $\dsm
\frac{\partial \varphi}{\partial \gamma}
=\frac{\sin\theta \cos\theta}
{\cos^2 \theta+\gamma^2 \sin^2 \theta}
$ 
\, thus: \,
$
\dsm \Absole{\frac{\partial \varphi}{\partial \gamma}}
\, \leq \,
\frac{1}{2\,\gamma}
\, \leq \,
1.2078
$
\, , \, and as before:
\[
-\frac{1}{\beta^2}\,\frac{\partial \beta}{\partial \gamma}
\,=\,
B'(\gamma)+
\nu\,  f_{\,3}'(\theta_\perp+s\,\nu)\,
\frac{\partial \varphi}{\partial \gamma}
\,=\,
B'(\gamma)\,+(\leq)\, 
\frac{(1-\gamma^2)\,M_{\,1}\,\delta}{2\,\gamma}
\,=\,
B'(\gamma)\,+(\leq)\, 0.08315
\]
where $0<s<1$. 
But for all $\gamma \in [\,\, 0.414\vge 1\,]$ we have:
\[
B'(\gamma)
\,=\, \int_0^{\frac{\pi}{2}}
\frac{\gamma\,\sin^2 t\,\,dt}
{\sqrt{\,\cos^2 t +\gamma^2\,\sin^2 t\,}^{\,3}\,} 
\,=\, 
\gamma\, F(\gamma)
\]
where $F$ decreases with $\gamma$\,, 
hence: \, $\absol{B'(\gamma)}\,\leq\, F(0.414)\,\gamma
\,\leq\, 5.3502$ \, , \, thus:
\[
\Absole{\frac{\partial \beta}{\partial \gamma}}
\,\leq\,
5.4334\,\beta_{max}^2
\,=\,
M_{\beta_\gamma}
\,\leq\,
2.317
\quad .
\]
If we set: \, $ \dsm
\overline{ E_\gamma}
\,=\,
\frac{\partial \beta}{\partial \gamma}\, 
\Big(\, \frac{\sin (2\,\theta) }
{\sqrt{\cos^2  \theta +\gamma^2\, \sin^2  \theta}}
\,+\,\int_{\theta}^{\theta_\perp}\frac{\cos(2\, t)\,dt}
{\sqrt{\,\cos^2 t+\gamma^2\,\sin^2 t\,}}
\,\Big)
-\frac{2}{\big( 1+\gamma \big)^2} 
$
\[
-\gamma\,\beta\,\Big(
\frac{\sin (2\,\theta)\, \sin^2  \theta }
{\sqrt{\cos^2  \theta +\gamma^2\, \sin^2  \theta}^{\,\,3}}
\,+\,
\int_{\theta}^{\theta_\perp}\frac{\cos(2\, t)\,\sin^2 t\,dt}
{\sqrt{\,\cos^2 t+\gamma^2\,\sin^2 t\,}^{\,\,3}}
\Big)
\,+\,
\frac{\beta\,\cos\big(2\, \theta_\perp\big)}
{\sqrt{\cos^2 \theta_\perp+\gamma^2\,\sin^2\theta_\perp }}\,
\frac{\partial \varphi}{\partial \gamma}
\]
we infer: \, $\dsm
\Absole{\frac{\partial E}{\partial \gamma}-\overline{ E_\gamma}}
\,\leq\,
M_{\beta_\gamma}\,\frac{\delta}{\gamma}
+\frac{\beta_{max}\,\delta}{\gamma^2}
\,+\,\frac{\beta_{max}\,\Delta_{y(1)}}{2\,\gamma}
\,\leq\,
0.3048
$ \, . \, We have moreover:
\[
\frac{\partial \beta}{\partial \gamma}
\,=\,
- \Big( B'(\gamma)\,+(\leq)\, 
\frac{(1-\gamma^2)\,M_{\,1}\,\delta}{2\,\gamma} \,\,\Big)
\times \Big(\,
\frac{1}{B(\gamma)}\,+(\leq)\,\Delta_\beta \Big)^2
\]
hence:
\[
\frac{\partial \beta}{\partial \gamma}
\,=\,
\widetilde{\beta}\,' (\gamma)   \,+(\leq)\, 
\frac{(1-\gamma^2)\,M_{\,1}\,\delta\,\beta_{max}^2}{2\,\gamma}
+ 5.4334 \cdot  2  \, \beta_{max} \, \Delta_\beta
\,=\,
\widetilde{\beta}\,' (\gamma)   \,+(\leq)\, \Delta_{\beta_\gamma} 
\]
where \, $\Delta_{\beta_\gamma} \leq 0.1519$ \, ,
thus if we set:
\[
\widetilde{ E_\gamma} (\gamma \vge \theta)
\,=\,
\widetilde{\beta}\, ' (\gamma)
\Big(\, \frac{\sin (2\,\theta) }
{\sqrt{\cos^2  \theta +\gamma^2\, \sin^2  \theta}}
\,+\,\int_{\theta}^{\theta_\perp}\frac{\cos(2\, t)\,dt}
{\sqrt{\,\cos^2 t+\gamma^2\,\sin^2 t\,}}
\,\Big)
-\frac{2}{\big( 1+\gamma \big)^2} 
\]
\[
-\gamma\,\widetilde{\beta}\,\Big(
\frac{\sin (2\,\theta)\, \sin^2  \theta }
{\sqrt{\cos^2  \theta +\gamma^2\, \sin^2  \theta}^{\,\,3}}
\,+\,
\int_{\theta}^{\theta_\perp}\frac{\cos(2\, t)\,\sin^2 t\,dt}
{\sqrt{\,\cos^2 t+\gamma^2\,\sin^2 t\,}^{\,\,3}}
\Big)
\,+\,
\frac{\widetilde{\beta}\,\cos\big(2\, \theta_\perp\big)}
{\sqrt{\cos^2 \theta_\perp+\gamma^2\,\sin^2\theta_\perp }}\,
\frac{\partial \varphi}{\partial \gamma}
\]
we get:  \, $\dsm
\Absole{\overline{ E_\gamma}-\widetilde{E_\gamma}}
\,\leq\,
\Delta_{ \beta_\gamma}\,  \frac{1+\pi}{\gamma}
+\Delta_{ \beta}\, \Big(\,\frac{1+\pi}{\gamma^2}+\frac{1}{2\,\gamma^2}\, \Big)
\,\leq\,
1.9631
$ \,
and finally:
\[
\Absole{\frac{\partial E}{\partial \gamma}(\gamma \vge \theta\vge \nu)
-\widetilde{E_\gamma}(\gamma \vge \theta)}
\,\leq\,
2.2678 \quad .
\]

\medskip

Now we estimate the partial derivatives of these $3$ functions of $2$
variables, starting by:
\[
\widetilde{E_\nu} (\gamma \vge \theta)
=
\frac{-\big(\,\widetilde{\beta}(\gamma)\,\big)^{2}}
{\sqrt{\cos^2\theta_\perp +\gamma^2\,\sin^2 \theta_\perp}}\,
\Big(\frac{\sin (2\,\theta) }
{\sqrt{\cos^2  \theta +\gamma^2\, \sin^2  \theta}}
+\int_{\theta}^{\theta_\perp}\frac{\cos(2\, t)\,dt}
{\sqrt{\,\cos^2 t+\gamma^2\,\sin^2 t\,}}\Big)
\]
\[
+\quad \frac{\widetilde{\beta}(\gamma)\,\cos(2\,\theta_\perp )}
{\sqrt{\,\cos^2 \theta_\perp +\gamma^2\,\sin^2 \theta_\perp}}
\]
where  \, $\widetilde{\beta}(\gamma)=1/B(\gamma)\, \leq \, 2/\pi$ \, 
and \, $\theta_\perp=\varphi(\gamma \vge \theta)$ \, , \, 
and we obtained:
\[
\Absole{B'(\gamma)}
\, \leq \,
\frac{\pi}{2\,\gamma^2}
\quad , \quad
\Absole{\frac{\partial \varphi}{\partial \theta}}
\, \leq \,
\frac{1}{\gamma}
\quad \hbox{and} \quad
\Absole{\frac{\partial \varphi}{\partial \gamma}}
\, \leq \,
\frac{1}{2\,\gamma}
\quad ,
\]
thus we get:
\quad $ \dsm
\Absole{\frac{\partial \widetilde{E_\nu} }{\partial \gamma}}
\, \leq \,
\Big(\,\frac{8}{\pi^2\,\gamma^3}+\frac{4}{\pi^2\,\gamma^2}+
\frac{1-\gamma^2}{\pi^2\,\gamma^4}\Big)\,\frac{\pi+1}{\gamma}
+\frac{4}{\pi^2\,\gamma}\,\Big(\frac{1}{\gamma^2}+
\frac{\pi+1}{\gamma^2}+\frac{1}{2\,\gamma^2}\Big)
$
\[
+\frac{2}{\pi\,\gamma^3}+\frac{2}{\pi\,\gamma^2}
+\frac{3+\gamma^2}{2\,\pi\,\gamma^4}
\quad \leq \quad
231
\]
and:
\quad $ \dsm
\Absole{\frac{\partial \widetilde{E_\nu}}{\partial \theta}}
\, \leq \,
\frac{2\,(1-\gamma^2)\,(\pi+1)}{\pi^2\,\gamma^5}
+\frac{4}{\pi^2\,\gamma}\,\Big(\,
\frac{1-\gamma^2}{2\,\gamma^3}
+\frac{1}{\gamma^2}+\frac{1}{\gamma}\,\Big)
+\frac{3+\gamma^2}{\pi\,\gamma^4}
\, \leq \,
107$ \quad .
\medskip

Then we have:
\[
\widetilde{E_\theta}
\, = \,
\frac{\widetilde{\beta}(\gamma)\,(1-\gamma^2)\,\sin^2(2\,\theta )} 
{2\,\sqrt{\cos^2  \theta +\gamma^2\, \sin^2  \theta}^{\,\,3}}
\, + \,
\frac{\widetilde{\beta}(\gamma)\,\cos(2\,\theta_\perp )}
{\sqrt{\,\cos^2 \theta_\perp+\gamma^2\,\sin^2 \theta_\perp}} \, \,
\frac{\partial \varphi}{\partial \theta}
\, - \,
\frac{\widetilde{\beta}(\gamma)\,\cos(2\,\theta )}
{\sqrt{\,\cos^2 \theta+\gamma^2\,\sin^2 \theta}} 
\]
where:
\, $\dsm
\frac{\partial \varphi}{\partial \theta}
=\frac{\gamma }{\cos^2 \theta + \gamma^2 \,\sin^2 \theta} 
$ \, , \, thus we get: 
$\dsm
\Absole{\frac{\partial^2 \varphi}{\partial \gamma\,\partial \theta}}
\,=\,
\Absole{\frac{\cos^2 \theta - \gamma^2 \,\sin^2 \theta }
{\big(\cos^2 \theta + \gamma^2 \,\sin^2 \theta\big)^2}  }
\,\leq\,
\frac{1}{\gamma^2}
$
by studying this function of $\theta$\,.
It comes therefore:
\[
\Absole{\frac{\partial \widetilde{E_\theta} }{\partial \gamma}}
\, \leq \,
\frac{1-\gamma^2}{\pi\,\gamma^5}
+\frac{2}{\pi\,\gamma^2}
+\frac{3\,(1-\gamma^2)}{\pi\,\gamma^4}
+\frac{2}{\pi\,\gamma^3}
+\frac{3+\gamma^2}{\pi\,\gamma^5}
+\frac{2}{\pi\,\gamma^3}
+\frac{2}{\pi\,\gamma^3}
+\frac{2}{\pi\,\gamma^2}
\, \leq \,
168 
\quad ,
\]
we also get:
\, $\dsm
\Absole{\frac{\partial^2 \varphi}{\partial \, \theta^2}}
\,=\,
\Absole{\frac{\gamma(1- \gamma^2)\,\sin (2\,\theta) }
{\big(\cos^2 \theta + \gamma^2 \,\sin^2 \theta\big)^2}  }
\,\leq\,
\frac{1- \gamma^2}{\gamma^3}
$ \, , \, and we infer:
\[
\Absole{\frac{\partial \widetilde{E_\theta} }{\partial \theta}}
\, \leq \,
\frac{4\,(1-\gamma^2)}{\pi\,\gamma^5} 
+\frac{3+\gamma^2}{\pi\,\gamma^4}
+\frac{2\,(1-\gamma^2)}{\pi\,\gamma^4}
+\frac{3+\gamma^2}{\pi\,\gamma^3}
\, \leq \,
156
\quad .
\]
\medskip

Finally, we have:
\[
\widetilde{ E_\gamma} (\gamma \vge \theta)
\,=\,
\widetilde{\beta}\, ' (\gamma)
\Big(\, \frac{\sin (2\,\theta) }
{\sqrt{\cos^2  \theta +\gamma^2\, \sin^2  \theta}}
\,+\,\int_{\theta}^{\theta_\perp}\frac{\cos(2\, t)\,dt}
{\sqrt{\,\cos^2 t+\gamma^2\,\sin^2 t\,}}
\,\Big)
-\frac{2}{\big( 1+\gamma \big)^2} 
\]
\[
-\gamma\,\widetilde{\beta}\,\Big(
\frac{\sin (2\,\theta)\, \sin^2  \theta }
{\sqrt{\cos^2  \theta +\gamma^2\, \sin^2  \theta}^{\,\,3}}
\,+\,
\int_{\theta}^{\theta_\perp}\frac{\cos(2\, t)\,\sin^2 t\,dt}
{\sqrt{\,\cos^2 t+\gamma^2\,\sin^2 t\,}^{\,\,3}}
\Big)
\,+\,
\frac{\widetilde{\beta}\,\cos\big(2\, \theta_\perp\big)}
{\sqrt{\cos^2 \theta_\perp+\gamma^2\,\sin^2\theta_\perp }}\,
\frac{\partial \varphi}{\partial \gamma}
\]
and here again the numeric approach will be useful to get a sharp enough estimate.
First, we get as above:
$\dsm \Absole{\widetilde{\beta}\,'(\gamma)}
\, \leq \,
\frac{2}{\pi\,\gamma^2}$
and
$\dsm \Absole{B''(\gamma)}
\, \leq \,
\frac{2\,\pi}{2\,\gamma^3}$\,,
hence
$\dsm \Absole{\widetilde{\beta}\,''(\gamma)}
\, \leq \,
\frac{8}{\pi^2\,\gamma^3}+\frac{4}{\pi\,\gamma^4}$\,\,.
Then we write patiently:
\[
\frac{\widetilde{\partial E_\gamma}}{\partial \gamma}
\,=\,
\frac{4}{\big( 1+\gamma \big)^2} 
+\widetilde{\beta}\, '' (\gamma)
\Big(\, \frac{\sin (2\,\theta) }
{\sqrt{\cos^2  \theta +\gamma^2\, \sin^2  \theta}}
\,+\,\int_{\theta}^{\theta_\perp}\frac{\cos(2\, t)\,dt}
{\sqrt{\,\cos^2 t+\gamma^2\,\sin^2 t\,}}
\,\Big)
\]
\[
-2\,\gamma\,\widetilde{\beta}\,' \,\Big(
\frac{\sin (2\,\theta)\, \sin^2  \theta }
{\sqrt{\cos^2  \theta +\gamma^2\, \sin^2  \theta}^{\,\,3}}
\,+\,
\int_{\theta}^{\theta_\perp}\frac{\cos(2\, t)\,\sin^2 t\,dt}
{\sqrt{\,\cos^2 t+\gamma^2\,\sin^2 t\,}^{\,\,3}}
\Big)
\]
\[
+\widetilde{\beta} \,\Big(
\frac{3\,\gamma^2\,\sin (2\,\theta)\, \sin^4  \theta }
{\sqrt{\cos^2  \theta +\gamma^2\, \sin^2  \theta}^{\,\,5}}
-\frac{\sin (2\,\theta)\, \sin^2  \theta }
{\sqrt{\cos^2  \theta +\gamma^2\, \sin^2  \theta}^{\,\,3}}
\Big)
\]
\[
+\widetilde{\beta} \,\Big(
\int_{\theta}^{\theta_\perp}\frac{3\,\gamma^2\,\cos(2\, t)\,\sin^4 t\,dt}
{\sqrt{\,\cos^2 t+\gamma^2\,\sin^2 t\,}^{\,\,5}}
-
\int_{\theta}^{\theta_\perp}\frac{\cos(2\, t)\,\sin^2 t\,dt}
{\sqrt{\,\cos^2 t+\gamma^2\,\sin^2 t\,}^{\,\,3}}
\Big)
\]
\[
\,-\,
\gamma\,\widetilde{\beta}\,
\frac{\cos\big(2\, \theta_\perp\big)\,\sin^2 \theta_\perp}
{\sqrt{\cos^2 \theta_\perp+\gamma^2\,\sin^2\theta_\perp }^3}\,
\frac{\partial \varphi}{\partial \gamma}
\,+\,
\frac{2\,\widetilde{\beta}\,'\,\cos\big(2\, \theta_\perp\big)}
{\sqrt{\cos^2 \theta_\perp+\gamma^2\,\sin^2\theta_\perp }}\,
\frac{\partial \varphi}{\partial \gamma}
\]
\[
\,-\,
\frac{\widetilde{\beta}\,\sin \big(2\, \theta_\perp\big)}
{\sqrt{\cos^2 \theta_\perp+\gamma^2\,\sin^2\theta_\perp }}\,
\Big(1+\gamma^2+\frac{1}{2}\,(1-\gamma^2)\,\cos \big( 2\,\theta_\perp\big)\,\Big)\,
\Big(\frac{\partial \varphi}{\partial \gamma}\Big)^2
\]
\[
\,+\,
\frac{\widetilde{\beta}\,\cos\big(2\, \theta_\perp\big)}
{\sqrt{\cos^2 \theta_\perp+\gamma^2\,\sin^2\theta_\perp }}\,
\frac{\partial^2 \varphi}{\partial \gamma^2}
\]
where:
\, $\dsm
\frac{\partial \varphi}{\partial \gamma}
=\frac{\sin\theta \cos\theta}
{\cos^2 \theta+\gamma^2 \sin^2 \theta}
$ 
\, and: \,
$
\dsm \Absole{\frac{\partial \varphi}{\partial \gamma}}
\, \leq \,
\frac{1}{2\,\gamma}
$
\, , \, 
thus:
$\dsm
\frac{\partial^2 \varphi}{\partial \gamma^2}
=\frac{-2\,\gamma\,\sin^3\theta \cos\theta}
{(\cos^2 \theta+\gamma^2 \sin^2 \theta)^2}
$\,.
The above estimates give: 
\[\frac{\widetilde{\partial E_\gamma}}{\partial \gamma}
\,=(\leq)\,
\frac{4}{\big( 1+\gamma \big)^2} 
+\Big(\,\frac{8}{\pi^2\,\gamma^3}+\frac{4}{\pi\,\gamma^4}\,\Big)
\Big(\, M_0
\,+\,\int_{0}^{\pi}\frac{\absol{\cos(2\, t)}\,dt}
{\sqrt{\,\cos^2 t+\gamma^2\,\sin^2 t\,}}
\,\Big)
\]
\[
+\frac{4}{\pi\gamma}\,\Big(\, M_5
\,+\,\int_{0}^{\pi}
\frac{\absol{\cos(2\, t)}\,\sin^2 t\,dt}
{\sqrt{\,\cos^2 t+\gamma^2\,\sin^2 t\,}^{\,\,3}}
\Big)
+\frac{2}{\pi}\, (M_5+M_6)
\]
\[
+\frac{2}{\pi}\,\Big(
\int_{0}^{\pi}\frac{3\,\absol{\cos(2\, t)}\,\sin^4 t\,dt}
{\sqrt{\,\cos^2 t+\gamma^2\,\sin^2 t\,}^{\,\,5}}
+
\int_{0}^{\pi}\frac{\absol{\cos(2\, t)}\,\sin^2 t\,dt}
{\sqrt{\,\cos^2 t+\gamma^2\,\sin^2 t\,}^{\,\,3}}
\Big)
+\frac{1}{\pi}\,\,M_7
\,+\,
\frac{M_4}
{\pi\,\gamma^4}
\]
\[
+\frac{3\,M_0}{4\,\pi\,\gamma^2}
\,+\,
\frac{4}{\pi}\,M_8
\]
where $M_8$ is the maximum of
$f_8 \dpe \dsm \theta \mapsto\frac{\absol{\sin^3\theta \cos\theta}}
{(\cos^2 \theta+\gamma^2 \sin^2 \theta)^2}
$
and similarly:
\[
f_4 \dpe \theta \mapsto
\frac{\absol{\cos\big(2\, \theta\big)}}
{\sqrt{\cos^2 \theta+\gamma^2\,\sin^2\theta }}
\quad , \quad
f_7 \dpe \theta \mapsto
\frac{\absol{\cos\big(2\, \theta\big)}\,\sin^2 \theta}
{\sqrt{\cos^2 \theta+\gamma^2\,\sin^2\theta }^3}\,
\]
\[
f_6 \dpe \theta \mapsto
\frac{3\,\absol{\sin (2\,\theta)}\, \sin^4  \theta }
{\sqrt{\cos^2  \theta +\gamma^2\, \sin^2  \theta}^{\,\,5}}
\quad , \quad
f_5 \dpe \theta \mapsto
\frac{\absol{\sin (2\,\theta)}\, \sin^2  \theta }
{\sqrt{\cos^2  \theta +\gamma^2\, \sin^2  \theta}^{\,\,3}}
\]
whose derivatives are bounded by $1/\gamma^6 < 200$.
The procedure $\star_{Maple}$ above answers:
$M_4 \leq 2.419$\,,
$M_5 \leq 4.323$\,,
$M_6 \leq 75.277$\,,
$M_7 \leq 14.096$\,,
$M_8 \leq 0.10862$ and we already knew that
$M_0 \leq 1.4145$\,, so we obtain finally:
\[
\Absole{\frac{\partial \widetilde{E_\gamma} }{\partial \gamma}}
\, \leq \,
453 \quad.
\]
At last, we have:
\[
\frac{\partial\widetilde{ E_\gamma} }{\partial \theta}
=
\widetilde{\beta}\, ' (\gamma)
\Big(\, \frac{(1-\gamma^2)\cos^2 \theta \, (2\,\cos^2 \theta-\gamma^2\,\sin^2 \theta) }
{\sqrt{\cos^2  \theta +\gamma^2\, \sin^2  \theta}^{\,3}}
\]
\[
+\frac{\cos(2\, \theta_\perp)}
{\sqrt{\,\cos^2 \theta_\perp+\gamma^2\,\sin^2 \theta_\perp\,}}\frac{\partial \varphi}{\partial \theta}
-\frac{\cos(2\, \theta)}
{\sqrt{\,\cos^2 \theta+\gamma^2\,\sin^2 \theta\,}}\,\Big)
\]
\[
-\gamma\,\widetilde{\beta}\,\Big(
\frac{3(1-\gamma^2)\sin^2 (2\,\theta)\, \sin^2  \theta }
{2\sqrt{\cos^2  \theta +\gamma^2\, \sin^2  \theta}^{\,\,5}}
+\frac{\cos(2\, \theta_\perp)\sin^2\theta_\perp}
{\sqrt{\,\cos^2 \theta_\perp+\gamma^2\,\sin^2 \theta_\perp\,}}
\frac{\partial \varphi}{\partial \theta}
-\frac{\cos(2\, \theta)\sin^2\theta}
{\sqrt{\,\cos^2 \theta+\gamma^2\,\sin^2 \theta\,}}\,\Big)
\]
\[
\,+\,
\frac{\widetilde{\beta}\,\sin\big(2\, \theta_\perp \big)\, (2\,\cos^2 \theta_\perp
-\gamma^2\sin^2 \theta_\perp)}
{\sqrt{\cos^2 \theta_\perp+\gamma^2\,\sin^2\theta_\perp }^{\,\,3}}\,
\,\frac{\partial \varphi}{\partial \gamma}\,\,\frac{\partial \varphi}{\partial \theta}\,
\,+\,
\frac{\widetilde{\beta}\,\cos\big(2\, \theta_\perp\big)}
{\sqrt{\cos^2 \theta_\perp+\gamma^2\,\sin^2\theta_\perp }}\,
\frac{\partial^2 \varphi}{\partial \gamma\partial \theta}
\]
and the rough estimates produce:
\[
\frac{\partial\widetilde{ E_\gamma} }{\partial \theta}
\,=(\leq)\,
\frac{2}{\pi\gamma^2}
\Big(\, \frac{1}{\gamma^{\,3}}+\frac{1}{\gamma^{\,2}}+\frac{1}{\gamma}\,\Big)
+\frac{2}{\pi}
\Big(\, \frac{1}{\gamma^{\,4}}+\frac{1}{\gamma}+1
+\frac{1}{2\,\gamma^{\,5}}+\frac{1}{\gamma^{\,3}}
\Big)
\,=(\leq)\, 144 \, \, ,
\]
which is good enough for us.
The Maple procedure from appendix \ref{approxderE}  
gives in 20 minutes the estimates used above:
\[
\Absole{\frac{\partial E}{\partial \nu}(\gamma \vge \theta \vge \nu) }
\, \leq \,
2.48
 \quad , \quad 
\Absole{\frac{\partial E}{\partial \theta}(\gamma \vge \theta \vge \nu) }
\, \leq \,
4.41
\quad \hbox{and} \quad
\Absole{\frac{\partial E}{\partial \gamma}(\gamma \vge \theta \vge \nu) }
\, \leq \,
4.33
\]
for all 
$(\gamma \vge \theta \vge \nu)\in [\, 0.414 \vge 1 \,]
\times[ \,0\vge \pi \,]\times [ \,-0.0354 \vge 0.0354 \,]$\,.

\section{Estimates on the partial derivatives of $cb$ and $en$} \label{der-iteration}

As in  section \ref{iteration}, we set  $h=1/(2\,s+1)$, 
$A=1+2\,h^2\,v$,
$B=1-2\,h^2\,v$
and
$C=1-2\,h^2\,u$\,, where $(u\vge v)\in D$ defined by: $\dsm
1 \,\leq \, u \,\leq \,
5
$
\, , \,
$\dsm
0 \,\leq \, v \,\leq \, u-1 \,\leq \, 4
$
\,
and if \, $u\geq 9/\,4$ \, :
\[
v 
\,\geq \, 
\sqrt{\frac{128\,u^2-144\,u-81-27\sqrt{9+32\,u}}{128}\,\,}
\quad ,
\]
as well as: \, $ \dsm
\alpha (u \vge v)=\frac{1}{\sqrt{u-v}}
$ 
\, and
\, $ \dsm
\beta  (u \vge v)=\frac{1}{\sqrt{u+v}}
$ \, .
The induction relation also writes
\quad $A\,\cos \theta_k\,\cos\theta_{k+1}
+
B\,\sin \theta_k\,\sin\theta_{k+1}
=C$ \, , and we set for each $\theta \in \MR$\,: \,
 $A\,\cos \theta\,\cos\theta^\dagger
+
B\,\sin \theta\,\sin\theta^\dagger
=0$ \,
to get as in section \ref{cb}:
\[
\theta^\dagger
=\theta+\frac{\pi}{2}-\arctan \Big( 
\frac{(A -B)\, \sin \theta\,  \cos \theta}
{A\,\cos^2 \theta + B \,\sin^2 \theta} \Big) 
=\theta+\frac{\pi}{2}-\arctan \Big( 
\frac{2\,h^2 v\, \sin 2 \theta}
{1+2\,h^2 v\, \cos 2 \theta} \Big) 
\]
and:
\[
\theta_{k+1} 
= \theta_k^\dagger 
-\arcsin \Big( \frac{C}{\sqrt{A^2\,\cos^2 \theta+B^2\sin^2 \theta}}\Big)
= \theta_k^\dagger 
-\arcsin \Big( \frac{1-2\,h^2 u}
{\sqrt{1+4\,h^4v^2+4\,h^2 v\,\cos 2 \theta}}\Big)
\]
hence: \quad $\theta_{k+1}=F (u\vge v\vge \theta_k)$ \quad 
for all $1\leq k \leq s$\,, where:
\[
F (u\vge v\vge \theta)=
\theta+\frac{\pi}{2}-\arctan \Big( 
\frac{2\,h^2 v\, \sin 2 \theta}
{1+2\,h^2 v\, \cos 2 \theta} \Big) 
-\arcsin \Big(\frac{1-2\,h^2 u}
{\sqrt{\rule{0cm}{0.4cm} 1+4\,h^4v^2+4\,h^2 v\,\cos 2 \theta}}\,\Big)
\, .
\]
It comes: \quad $ \dsm
\frac{\partial F}{\partial u}
\,=\,
\frac{h}
{\sqrt{u+v\,\cos 2\theta -h^2\,\big(u^2-v^2\big)}}
\,=(\leq)\, m_u\, h = M_u
$\quad ,
\[
\frac{\partial F}{\partial v}
\,=\,
\frac{2\,h^3 v+h\, \cos 2 \theta}
{(1+4\,h^4v^2+4\,h^2 v\,\cos 2 \theta)\,
\sqrt{u+v\,\cos \big(2\,\theta\big)-h^2\,(u^2-v^2)}}
\]
\[
-\frac{2\,h^2 \sin 2 \theta}
{1+4\,h^2 v\, \cos 2 \theta+4\,h^4 v^2}
\,=\, h \, F_v(u\vge v \vge \theta)\,=(\leq)\, m_v\, h = M_v
\quad \hbox{and:}
\]
\[
\frac{\partial F}{\partial \theta}
\,=\,
\frac{-2\,h\, v\,\sin 2 \theta}
{(1+4\,h^4v^2+4\,h^2 v\,\cos 2 \theta)\,
\sqrt{u+v\,\cos \big(2\,\theta\big)-h^2\,(u^2-v^2)}}
\]
\[
+1-\frac{4\,h^2 v\,\cos 2 \theta+8\,h^4 v^2 }
{1+4\,h^2 v\, \cos 2 \theta+4\,h^4 v^2}
\,=\, 1+h \, F_\theta(u\vge v \vge \theta)\,=\, 1 \,+(\leq)\,m_\theta h
=M_\theta
\]
where \quad $\dsm m_u=\frac{1}{\sqrt{(u-v)\big(1-h^2 (u+v)\big)}}
\leq \frac{1}{\sqrt{\rule{0cm}{0.4cm} 1-9\,h^2}}$\quad ,
but we need sharp estimates on $m_v$ and $m_\theta$ and we 
will again obtain them numerically. The induction relation implies:
\[
\frac{\partial \theta_{k+1}}{\partial \theta}=
\frac{\partial F}{\partial \theta}
(u\vge v \vge \theta_k) \,
\frac{\partial \theta_{k}}{\partial \theta}
\quad ,
\]
\[
\frac{\partial \theta_{k+1}}{\partial u}=
\frac{\partial F}{\partial u}(u\vge v \vge \theta_k)
+\frac{\partial F}{\partial \theta} (u\vge v \vge \theta_k)\,
\frac{\partial \theta_{k}}{\partial u}
\]
and:
\[
\frac{\partial \theta_{k+1}}{\partial v}=
\frac{\partial F}{\partial v}(u\vge v \vge \theta_k)
+\frac{\partial F}{\partial \theta} (u\vge v \vge \theta_k)\,
\frac{\partial \theta_{k}}{\partial v}
\]
for all $1\leq k \leq s$\,, thus:
\[
\Absole{\frac{\partial \theta_k}{\partial \theta}}
\,\leq\,
M_\theta^{\,k-1}
\quad ,\quad 
\Absole{\frac{\partial \theta_{k}}{\partial u}}
\,\leq\,
M_u\,\frac{M_{\theta}^{\,k-1}-1}{M_{\theta}-1}
\quad \hbox{and} \quad 
\Absole{\frac{\partial \theta_{k}}{\partial v}}
\,\leq\,
M_v\,\frac{M_{\theta}^{\,k-1}-1}{M_{\theta}-1}
\quad , 
\]
hence if $k\geq 2$: \qquad \qquad \qquad \quad 
$\dsm
\Absole{\frac{\partial \theta_k}{\partial \theta}}
\,\leq\,
\Big(1+\frac{m_\theta}{2\,s+1}\Big)^{k-1}
$ \quad ,
\[
\Absole{\frac{\partial \theta_{k}}{\partial u}}
\,\leq\,
\frac{m_u}{m_{\theta}}\,\,
\Big(1+\frac{m_\theta}{2\,s+1}\Big)^{k-1}
\qquad \hbox{and} \qquad 
\Absole{\frac{\partial \theta_{k}}{\partial v}}
\,\leq\,
\frac{m_v}{m_{\theta}}\,\,
\Big(1+\frac{m_\theta}{2\,s+1}\Big)^{k-1}
\quad . 
\]
For all $(u\vge v \vge \theta) \in K=D \times [0\vge \pi ]$ we have:
\[
cb(u\vge v\vge \theta)=
\frac{1}{\sqrt{u-v}}\,   \cos \theta  \,  \cos \theta_{s+1}  + 
\frac{1}{\sqrt{u+v}}\,    \sin \theta  \, \sin \theta_{s+1}
-\frac{1}{2\,s+1}
\]
and:
\qquad $ en(u\vge v\vge \theta_1)\,\,=$
\[
\frac{1}{2\,s+1}
\sum_{k=1}^{s} \cos\big(\theta_{k+1}+\theta_k\big)
+\frac{1}{\sqrt{u+v}}\,\sin \theta \, \sin \theta_{s+1}
+
\frac{s\,\sqrt{u+v}-\big(s+1\big)\,\sqrt{u-v}}
{\big(\sqrt{u+v}+\sqrt{u-v}\big)\,\big(2s+1\big)}  
\]
hence:
\[
\Absole{\frac{\partial cb}{\partial \theta}}
\,\leq\,
1+\Absole{\frac{\partial \theta_{s+1}}{\partial \theta}}
\,\leq\,
1+\Big(1+\frac{m_\theta}{2\,s+1}\Big)^{s}
\,\leq\,
1+e^{\frac{m_\theta}{2}}
\]
and:
\[
\Absole{\frac{\partial en}{\partial \theta}}
\,\leq\,
1+\Absole{\frac{\partial \theta_{s+1}}{\partial \theta}}
+\frac{1}{2s+1}\,\sum_{k=1}^{s} 
\,\Big(\,\Absole{\frac{\partial \theta_{k} }{\partial \theta}}
+\Absole{\frac{\partial \theta_{k+1} }{\partial \theta}}\,\Big)
\]
thus:
\[
\Absole{\frac{\partial en}{\partial \theta}}
\,\leq\,
1+e^{\frac{m_\theta}{2}}
+\frac{2}{2s+1}\,\sum_{k=1}^{s+1} 
\,\,\Absole{\frac{\partial \theta_{k} }{\partial \theta}}
\,\leq\,
1+e^{\frac{m_\theta}{2}}
\,\,\Big(\,1+\frac{2}{m_{\theta}}\,e^{\frac{m_\theta}{2\,(2\,s+1)}}\,\Big)
\quad .
\]
Similarly, $u-v \geq 1$ yields:
\[
\Absole{\frac{\partial cb}{\partial u}}
\,\leq\,
\frac{1}{2}
+\frac{m_u}{m_{\theta}}\,\,
\Big(1+\frac{m_\theta}{2\,s+1}\Big)^{s}
\,\leq\,
\frac{1}{2}
+\frac{m_u}{m_{\theta}}\,\,e^{\frac{m_\theta}{2}}
\quad 
\]
and
\[
\Absole{\frac{\partial cb}{\partial v}}
\,\leq\,
\frac{1}{2}
+\frac{m_v}{m_{\theta}}\,\,
\Big(1+\frac{m_\theta}{2\,s+1}\Big)^{s}
\,\leq\,
\frac{1}{2}
+\frac{m_v}{m_{\theta}}\,\,e^{\frac{m_\theta}{2}}
\quad ,
\]
then:
\[
\Absole{\frac{\partial en}{\partial u}}
\,\leq\,
\frac{2}{3}
+\frac{m_u}{m_{\theta}}\,\,
\Big(1+\frac{m_\theta}{2\,s+1}\Big)^{s}
+\frac{2}{2s+1}\,\sum_{k=1}^{s+1} 
\,\,\frac{m_u}{m_{\theta}}\,\,
\Big(1+\frac{m_\theta}{2\,s+1}\Big)^{k-1}
\]
\[
\,\leq\,
\frac{2}{3}
+\frac{m_u}{m_{\theta}}\,\,e^{\frac{m_\theta}{2}}
\,\,\Big(\,1+\frac{2}{m_{\theta}}\,e^{\frac{m_\theta}{2\,(2\,s+1)}}\,\Big)
\]
and
\[
\Absole{\frac{\partial en}{\partial v}}
\,\leq\,
\frac{7}{4}
+\frac{m_v}{m_{\theta}}\,\,
\Big(1+\frac{m_\theta}{2\,s+1}\Big)^{s}
+\frac{2}{2s+1}\,\sum_{k=1}^{s+1} 
\,\,\frac{m_v}{m_{\theta}}\,\,
\Big(1+\frac{m_\theta}{2\,s+1}\Big)^{k-1}
\]
\[
\,\leq\,
\frac{7}{4}
+\frac{m_v}{m_{\theta}}\,\,e^{\frac{m_\theta}{2}}
\,\,\Big(\,1+\frac{2}{m_{\theta}}\,e^{\frac{m_\theta}{2\,(2\,s+1)}}\,\Big)
\quad .
\]

The steady computation of the second derivatives of $f$ gives moreover:
\[
\Absole{\frac{\partial^2 F}{\partial v \partial u}}
\,\leq\,
\frac{(1+8\,h^2)}{\big(1-9\,h^2\big)^{3/2}}\,\,h
\quad ,\quad 
\Absole{\frac{\partial^2 F}{\partial \theta \partial u}}
\,\leq\,
\frac{8}{\big(1-9\,h^2\big)^{3/2}}\,\,h
\quad ,
\]
\[
\Absole{\frac{\partial^2 F}{\partial v^2}}
\,\leq\,
\frac{4\,h^2\,\big(1-2\,h^2\big)\,\big(1+8\,h^2\big)^2}
{\big(1-8\,h^2\big)^4\,\sqrt{1-9\,h^2}}\,\,h
+
\frac{\big(1-2\,h^2\big)\,\big(1+8\,h^2\big)^2}
{2\,\big(1-8\,h^2\big)^2\,\big(1-9\,h^2\big)^{3/2}}\,\,h
\]
\[
+
\frac{2\,h^2\,\big(1-2\,h^2\big)}
{\big(1-8\,h^2\big)^2\,\sqrt{1-9\,h^2}}\,\,h
+
\frac{8\,h^3\,\big(1+2\,h^2\big)}{\big(1-8\,h^2\big)^4}\,\,h
\quad ,
\]
\[
\Absole{\frac{\partial^2 F}{\partial \theta \partial v}}
\,\leq\,
\frac{2\,\big(1-2\,h^2\big)}
{\big(1-8\,h^2\big)^2\,\sqrt{1-9\,h^2}}\,\,h
+
\frac{32\,h^2\big(1-2\,h^2\big)\,\big(1+8\,h^2\big)}
{\big(1-8\,h^2\big)^4\,\sqrt{1-9\,h^2}}\,\,h
\]
\[
+
\frac{8\,\big(1-2\,h^2\big)\,\big(1+8\,h^2\big)}
{\big(1-8\,h^2\big)^2\,\big(1-9\,h^2\big)^{3/2}}\,\,h
+
\frac{4\,h}{\big(1-8\,h^2\big)^2}\,\,h
+
\frac{64\,h^3}{\big(1-8\,h^2\big)^4}\,\,h
\]
and finally:
\quad $\dsm
\Absole{\frac{\partial^2 F}{\partial \theta^2}}
\,\leq\,
\frac{16\,\big(1-2\,h^2\big)}
{\big(1-8\,h^2\big)^2\,\sqrt{1-9\,h^2}}\,\,h
+\frac{1024\,h^2\big(1-2\,h^2\big)}
{\big(1-8\,h^2\big)^4\,\sqrt{1-9\,h^2}}\,\,h
$
\[
+\frac{64\,\big(1-2\,h^2\big)}
{\big(1-8\,h^2\big)^2\,\big(1-9\,h^2\big)^{3/2}}\,\,h
+
\frac{32\,h}{\big(1-8\,h^2\big)^2}\,\,h
+
\frac{128\,h^3\,\big(1+8\,h^2\big)}{\big(1-8\,h^2\big)^4}\,\,h
\quad ,
\]
and if $4\leq s \leq 14$ we get $h \leq 1/9$ hence:
\[
m_u\, = \,  \frac{1}
{\sqrt{\big(u-v\big)\,\big(1-h^2\,(u+v)\big)}}
\,\leq\, 
1.061
\]
and:
\qquad \quad 
$\dsm
\Absole{\frac{\partial^2 F}{\partial v \partial u}}
\,\leq\,
0.104
\quad ,\quad 
\Absole{\frac{\partial^2 F}{\partial \theta \partial u}}
\,\leq\,
0.963
\quad , \quad
\Absole{\frac{\partial^2 F}{\partial v^2}}
\,\leq\,
0.112
\quad ,
$
\[
\Absole{\frac{\partial^2 F}{\partial \theta \partial v}}
\,\leq\,
16.02
\quad \hbox{and} \quad
\Absole{\frac{\partial^2 F}{\partial \theta^2}}
\,\leq\,
15.18
\quad ,
\]
The C procedure from appendix \ref{mvmt} 
estimates $m_v$ and $m_\theta$ with a step of
$1/(4\,n)$\,, hence and uncertainty of $2.3/n$ on $m_v$ and $20.2/n$
on $m_\theta$\,. 
If $n=100$\,, it answers  $m_v=1.042$  and  $m_\theta=5.272$
in 40 minutes, hence we get:
\[
m_u
\,\leq \, 
1.061
\quad , \quad
m_v
\,\leq \, 
1.065
\quad \et \quad
m_\theta
\,\leq \, 
5.474
\quad ,
\]
which leads to the estimates used above:
\[
\Absole{\frac{\partial cb}{\partial \theta}}
\,\leq\,
16.5
\quad
,
\quad
\Absole{\frac{\partial cb}{\partial u}}
\,\leq\,
3.5
\quad
\hbox{and}
\quad
\Absole{\frac{\partial cb}{\partial v}}
\,\leq\,
3.6
\quad
,
\]
\[
\Absole{\frac{\partial en}{\partial \theta}}
\,\leq\,
24.1
\quad
,
\quad
\Absole{\frac{\partial en}{\partial u}}
\,\leq\,
5.15
\quad
\hbox{and}
\quad
\Absole{\frac{\partial en}{\partial v}}
\,\leq\,
6.25
\quad .
\]
If $s=3$\,, the same $n=100$ leads in $4$ minutes to:
\[
m_u
\,\leq \, 
1.11
\quad , \quad
m_v
\,\leq \, 
1.12
\quad \et \quad
m_\theta
\,\leq \, 
6.16
\quad ,
\]
hence:
\[
\Absole{\frac{\partial cb}{\partial \theta}}
\,\leq\,
22.2
\quad
,
\quad
\Absole{\frac{\partial cb}{\partial u}}
\,\leq\,
4.35
\quad
\hbox{and}
\quad
\Absole{\frac{\partial cb}{\partial v}}
\,\leq\,
4.4
\]
\[
\Absole{\frac{\partial en}{\partial \theta}}
\,\leq\,
32.9
\quad
,
\quad
\Absole{\frac{\partial en}{\partial u}}
\,\leq\,
6.46
\quad
\hbox{and}
\quad
\Absole{\frac{\partial en}{\partial v}}
\,\leq\,
6.9
\quad .
\]


If $s=2$\,, we have $h=0.2$ hence:
\[
\Absole{\frac{\partial^2 F}{\partial v \partial u}}
\,\leq\,
0.52
\quad ,\quad 
\Absole{\frac{\partial^2 F}{\partial \theta \partial u}}
\,\leq\,
3.12
\quad , \quad 
\Absole{\frac{\partial^2 F}{\partial v^2}}
\,\leq\,
1.38
\quad ,
\]
\[
\Absole{\frac{\partial^2 F}{\partial \theta \partial v}}
\,\leq\,
12.68
\quad\hbox{and} \quad
\Absole{\frac{\partial^2 F}{\partial \theta^2}}
\,\leq\,
100.33
\quad ,
\]
thus the uncertainty on $m_v$ equals 
$9.2/n$ and the one on $m_\theta$ equals $74.5/n$\,,  
so we choose $n=250$ to get the estimates
$m_v=1.161$ and $m_\theta=7.474$, hence we obtain:
\[
m_u \,= \, 1.25
\quad , \quad
m_v \, = \, 1.2
\quad \et \quad
m_\theta \, = \, 7.8
\quad ,
\]
but the above estimates on the derivatives of $en$ and $cb$ 
would require $n=500$ in section \ref{iteration}, hence 
$8$ hours of computation in order to conclude.
Fortunately, we can easily sharpen them: we have
$M_\theta = 2.56$\,,  $M_u=0.25$  and  $M_v=0.25$\,, as well as:
\[
\Absole{\frac{\partial \theta_2}{\partial \theta}}
\,\leq\,
M_\theta
\quad \et \quad 
\Absole{\frac{\partial \theta_3}{\partial \theta}}
\,\leq\,
M_\theta^{\,2}
\quad , 
\]
thus it comes directly:
\quad $ \dsm
\Absole{\frac{\partial cb}{\partial \theta}}
\,\leq\,
1+\Absole{\frac{\partial \theta_{3}}{\partial \theta}}
\,\leq\,
1+2.56^2
\,\leq\,
7.56
$ \quad
and:
\[
\Absole{\frac{\partial en}{\partial \theta}}
\,\leq\,
1+\Absole{\frac{\partial \theta_{3}}{\partial \theta}}
+\frac{1}{5}\,\Big(\,1
+2\,\Absole{\frac{\partial \theta_{2} }{\partial \theta}}
+\Absole{\frac{\partial \theta_{3} }{\partial \theta}}\,\Big)
\,\leq\,
10.09 \quad .
\]
Moreover, if $w\in \{u\vge v\}$ we have: 
\, $ \dsm
\Absole{\frac{\partial \theta_2}{\partial w}}
\,\leq\,
M_w
$
\, and \, $ \dsm
\Absole{\frac{\partial \theta_3}{\partial \theta}}
\,\leq\,
M_w\,(1+M_\theta)
$
 \, , \, hence:
\[
\Absole{\frac{\partial cb}{\partial u}}
\,\leq\,
\frac{1}{2}
+M_u\,(1+M_\theta)
\,\leq\,
1.39
\quad \et \quad 
\Absole{\frac{\partial cb}{\partial v}}
\,\leq\,
\frac{1}{2}
+M_v\,(1+M_\theta)
\,\leq\,
1.39 \quad ,
\]
then:
\[
\Absole{\frac{\partial en}{\partial u}}
\,\leq\,
\frac{2}{3}+
\Absole{\frac{\partial \theta_{3}}{\partial u}}
+\frac{1}{5}\,\Big(\,1
+2\,\Absole{\frac{\partial \theta_{2} }{\partial u}}
+\Absole{\frac{\partial \theta_{3} }{\partial u}}\,\Big)
\,\leq\,
2.22
\]
and
\[
\Absole{\frac{\partial en}{\partial v}}
\,\leq\,
\frac{7}{4}+
\Absole{\frac{\partial \theta_{3}}{\partial v}}
+\frac{1}{5}\,\Big(\,1
+2\,\Absole{\frac{\partial \theta_{2} }{\partial v}}
+\Absole{\frac{\partial \theta_{3} }{\partial v}}\,\Big)
\,\leq\,
2.81 \quad .
\]
These estimates used in section \ref{iteration} finish
our proof of Gr{\"u}nbaum conjecture.

\section{Appendix: procedures}

\subsection{Estimating $\mu (0.0354)$ in section \ref{numerique-ed}}
\label{procmu}

This Maple procedure:
\begin{verbatim}
delta:=0.0354; 
coefferr:=evalf(2.48*delta/2+4.41*Pi/(2*79)+4.33*0.586/(2*64));
n:=6; err:= evalf(coefferr/n);  mm:=1000:
for igamma from 0 to 64*n do gammaa:=evalf(0.414+0.586*igamma/(64*n)): 
for itheta from 0 to 158*n-1 do theta:=evalf(Pi*itheta/(158*n)):
thetaperp:=evalf(theta+Pi/2-arctan((1-gammaa)*sin(theta)*cos(theta)
/(1-(1-gammaa)*(sin(theta))^2))):
for inu from -n to n do y1:=evalf(thetaperp+delta*inu/n):
beta:=evalf(1/int(1/sqrt(1-(1-gammaa^2)*(sin(t))^2), t=theta..y1)):
alpha:=beta/gammaa: 
if alpha<=1 then if beta<=alpha then if alpha+beta>=4/3 then 
E:=evalf((1-gammaa)/(1+gammaa)+beta*(
sin(2*theta)/sqrt((1-(1-gammaa^2)*(sin(theta))^2))
+int(cos(2*t)/sqrt(1-(1-gammaa^2)*(sin(t))^2),t=theta..y1))):
if E < mm then mm:=E: fi:fi:fi:fi:od:od:od: mm; mu:=evalf(mm-err);
\end{verbatim} 
realizes an uncertainty of \, $0.1514/n$ \, 
on the minimum $\mu (\delta)$ of the function $E$ on the domain $\Delta$\,.

\subsection{Estimating the partial derivative of $E$
in section \ref{deriveesEN}}\label{approxderE}

This Maple procedure:
\begin{verbatim}
n:=500; Mg:=0: Mz:=0: Mnu:=0:
for ig from 0 to n do g:=evalf(0.414+0.586*ig/(n)): 
betatilde:=evalf(1/InverseJacobiAM(Pi/2,sqrt(1-g^2))):
for iz from 0 to 2*n-1 do z:=evalf(Pi*iz/(2*n)): 
zperp:=evalf(z+Pi/2-arctan((1-g)*0.5*sin(2*z)/(1-(1-g)*(sin(z))^2))):
II:=evalf(int(cos(2*t)/sqrt(1-(1-g^2)*(sin(t))^2),t=z..zperp)+
sin(2*z)/sqrt(1-(1-g^2)*(sin(z))^2)):
dnuen:=abs(evalf(-betatilde^2*II/sqrt(1-(1-g^2)*(sin(zperp))^2))
+betatilde*cos(2*zperp)/sqrt(1-(1-g^2)*(sin(zperp))^2)):
dphiz:=evalf(g/((cos(z))^2+g^2*(sin(z))^2)):
dzen:=abs(evalf( betatilde*cos(2*zperp)*dphiz/sqrt(1-(1-g^2)
*(sin(zperp))^2)
+betatilde*(1-g^2)*(sin(2*z))^2/(2*(sqrt(1-(1-g^2)*(sin(z))^2))^3)
-betatilde*cos(2*z)/sqrt(1-(1-g^2)*(sin(z))^2))):
dphig:=evalf(sin(z)*cos(z)/((cos(z))^2+g^2*(sin(z))^2)):
dbg:=evalf(betatilde^2*g
*int((sin(t))^2/(1-(1-g^2)*(sin(t))^2)^(3/2),t=z..zperp)
-betatilde^2*dphig/sqrt(1-(1-g^2)*(sin(zperp))^2)):
III:=evalf(int(cos(2*t)*(sin(t))^2/(1-(1-g^2)*(sin(t))^2)^(3/2),
 t=z..zperp)
+sin(2*z)*(sin(z))^2/(1-(1-g^2)*(sin(z))^2)^(3/2)):
dgen:=abs(evalf(dbg*II-2/(1+g)^2-g*betatilde*III
+betatilde*cos(2*zperp)*dphig/sqrt(1-(1-g^2)*(sin(zperp))^2))):
if dnuen >Mnu then Mnu:=dnuen:fi: if dzen >Mz then Mz:=dzen:fi:  
if dgen >Mg then Mg:=dgen:fi:     od:od:       Mnu;Mz;Mg;
MMnu:=evalf(Mnu+231*0.586/n+107*Pi/(2*n));
MMz:=evalf(Mz+231*0.586/n+168*Pi/(2*n));
MMg:=evalf(Mg+144*0.586/n+453*Pi/(2*n));
MderEnu:= evalf(MMnu+1.04753); 
MderEtheta:= evalf(MMz+1.53551);
MderEgamma:= evalf(MMg+2.26776);
\end{verbatim} 
estimates the maxima of $\widetilde{E}_\nu$\,,
$\widetilde{E}_\gamma$ and $\widetilde{E}_\theta$\,,
then the partial derivatives of the function $E$.

\subsection{Estimating the minimum of $m$ in section \ref{iteration}}
\label{minimumm}

This procedure in C estimates the minimum of
$m$ with a step of $1/(4\,n)$\,:
\begin{verbatim}
#include <math.h>    #include <stdio.h>     
int main (void)
{ int s; int n; double m; double dm; double h; int iu; double u; 
  double vmin; int nv; int iv; double v; double A; double B; double C; 
  double a; double b; int iz; double z;  double x; double y; double xx; 
  double yy; double en; double cb; int k; double mm; m = 0.1001e3;         
n = 121; 
for (s = 4; s <= 14 ; s++)
{h = (double)(0.1e1 / (double) (2 * s + 1));
  for (iu = 0; iu <= 5 * n; iu++)
  { u = (double)(0.1e1 + (double) (iu / (double) n) / (double) 0.4e1);
    C = 0.1e1 + (-0.2e1) * h * h * u;
  for (iv = 0; iv <= iu; iv++)
  { v = (double)((double) (iv / (double) n) / (double) 0.4e1);
    A = 0.1e1 + 0.2e1 * h * h * v; B = 0.1e1 + (-0.2e1) * h * h * v;
    a = (double)(0.1e1 / (double) sqrt(u - v)); 
    b = (double)(0.1e1 / (double) sqrt(u + v));
  for (iz = 0; iz <= 13 * n - 1; iz++)
  { z = (double)(0.3141592654e1 * (double) iz / 
        (double) n / (double) 0.13e2);
    x = (double)(A * cos(z)); y = (double)(B * sin(z));
    en = (double)(((double) s * a - (double) (s + 1) * b) / 
         (double) (2 * s + 1) / (double) (a + b));
  for (k = 1; k <= s; k++)
  { xx = (double) (A*(C*x-y*sqrt(x*x+y*y + (-0.1e1) * C * C)) 
         / (double) (x * x + y * y));
    yy = (double) (B*(C*y+x*sqrt(x*x+y*y + (-0.1e1) * C * C)) 
         / (double) (x * x + y * y));
    en = en + (double)((x * xx / (double) A / (double) A + 
         (-0.1e1) * y * yy / (double) B / (double) B) / 
         (double) (2 * s + 1));
    x = xx; y = yy; }
    en = (double)(en + b * sin(z) * y / (double) B);
    cb = (double)(a * cos(z) * x / (double) A + 
          b * sin(z) * y / (double) B - 0.1e1 / (double) (2 * s + 1));
    mm = (fabs(cb) >= fabs(en) ? fabs(cb) : fabs(en));
    if (mm < m) m = mm; } } }
  for (iu = 5 * n + 1; iu <= 16 * n; iu++)
  { u = (double)(0.1e1 + (double) (iu / (double) n) / (double) 0.4e1);
    C = 0.1e1 + (-0.2e1) * h * h * u;
    vmin = (double)(sqrt(0.1e1 * u * u - 0.9e1 / (double) 0.8e1 * u 
           - 0.81e2 / (double) 0.128e3 - 0.27e2 / 0.128e3 * 
           sqrt(0.9e1 + 0.32e2 * u)));
    nv = (double)(ceil(0.4e1 * (double) n * (u - 0.1e1 - vmin)) + 0.1e1);
  for (iv = 0; iv <= (int) nv; iv++)
  { v = (double)(vmin + (double) iv * (u - 0.1e1 - vmin) / (double) nv);
    A = 0.1e1 + 0.2e1 * h * h * v; B = 0.1e1 + (-0.2e1) * h * h * v;
    a = (double)(0.1e1 / (double) sqrt(u - v));
    b = (double)(0.1e1 / (double) sqrt(u + v));
  for (iz = 0; iz <= 13 * n - 1; iz++)
   { z = (double)(0.3141592654e1 * (double) iz / 
        (double) n / (double) 0.13e2);
    x = (double)(A * cos(z)); y = (double)(B * sin(z));
    en = (double)(((double) s * a - (double) (s + 1) * b) / 
         (double) (2 * s + 1) / (double) (a + b));
  for (k = 1; k <= s; k++)
  { xx = (double) (A*(C*x-y*sqrt(x*x+y*y + (-0.1e1) * C * C)) 
         / (double) (x * x + y * y));
    yy = (double) (B*(C*y+x*sqrt(x*x+y*y + (-0.1e1) * C * C)) 
         / (double) (x * x + y * y));
    en = en + (double)((x * xx / (double) A / (double) A + 
         (-0.1e1) * y * yy / (double) B / (double) B) / 
         (double) (2 * s + 1));
    x = xx; y = yy; }
    en = (double)(en + b * sin(z) * y / (double) B);
    cb = (double)(a * cos(z) * x / (double) A + 
          b * sin(z) * y / (double) B - 0.1e1 / (double) (2 * s + 1));
    mm = (fabs(cb) >= fabs(en) ? fabs(cb) : fabs(en));
    if (mm < m) m = mm; } } }
}
  printf("m = %.10f\n", m);
  cb = (double)( (0.165e2) / (double) (0.8e1) +(0.35e1) /(double) (0.8e1)
  +(double) (0.36e1) *(0.3141592654e1)/ (double) (0.26e2) )/ (double) n;
  en = (double)( (0.241e2)/ (double) (0.8e1) +(0.515e1) /(double) (0.8e1)
  +(double)  (0.625e1)*(0.3141592654e1)/ (double) 0.26e2)/ (double) n; 
  dm  = (fabs(cb) >= fabs(en) ? fabs(cb) : fabs(en));
  printf("dm = %.10f\n", dm);  printf("m-dm = %.10f\n", m-dm);
return(0);}
\end{verbatim} 
as well as the uncertainty $\delta m$\,.

\subsection{Estimating $m_v$ and $m_\theta$ in section \ref{der-iteration}}
\label{mvmt}

The following C procedure estimates $m_v$ and $m_\theta$ with a step of
$1/(4\,n)$\,, hence and uncertainty of $2.3/n$ on $m_v$ and $20.2/n$ on $m_\theta$\,:
\begin{verbatim}
#include <math.h>
#include <stdio.h>
int main (void)
{ int s; int n; double mv; double mx; double c;
  int iu; double u; double vmin; int nv; int iv; double v;
  int ix; double x; double dfv; double dfx;
n = 100; mv = 0; mx = 0;
for (s=4; s <= 14 ; s++)
{c = (double)(0.1e1 / (double) (2 * s + 1));
  for (iu = 0; iu <= 5 * n; iu++)
  {u = (double) (0.1e1 + (double) (iu / (double) n) / (double) 0.4e1);
  for (iv = 0; iv <= iu; iv++)
  {v = (double)((double) (iv / (double) n) / (double) 0.4e1);
  for (ix = 0; ix <= 13 * n - 1; ix++)
  {x = (double)(0.3141592654e1 *(double) ix /(double) n /(double)0.13e2);
  dfv = fabs((0.1e1 + (-0.2e1) * c * c * u) * (cos(0.2e1 * x) + 0.2e1 
  * c * c * v) / (double) (0.1e1 + 0.4e1 * c * c * v * cos(0.2e1 * x) 
  + 0.4e1 *  c * c * c * c * v * v) / (double) sqrt(u + v * cos(0.2e1 
  * x) + (-0.1e1) * c * c * u * u + c * c * v * v) - 0.2e1 * c * 
  sin(0.2e1 * x) / (double) (0.1e1 + 0.4e1 * c * c * v * 
  cos(0.2e1 * x) + 0.4e1 * c * c * c * c * v * v));
  dfx = fabs(0.2e1 * v * (0.1e1 + (-0.2e1) * c * c * u) * 
  sin(0.2e1 * x) / (0.1e1 + 0.4e1 * c * c * v * cos(0.2e1 * x) 
  + 0.4e1 * c * c * c * c * v * v) / (double) sqrt(u + v * cos(0.2e1 
  * x) + (-0.1e1) * c * c * u * u + c * c * v * v) +  0.4e1 * c * v * 
  (cos(0.2e1 * x) + 0.2e1 * c * c * v) / (double) (0.1e1 + 0.4e1 * c 
  * c * v * cos(0.2e1 * x) + 0.4e1 * c * c * c  * c * v * v));
     if ((double) mv < dfv) mv = (double) dfv;
     if ((double) mx < dfx) mx = (double) dfx; } } }
  for (iu = 5 * n + 1; iu <= 16 * n; iu++)
  {u = (double)(0.1e1 + (double) (iu / (double) n) / (double) 0.4e1);
   vmin = (double)(sqrt(0.1e1 * u * u - 0.9e1 / (double) 0.8e1 * u 
   - 0.81e2 / (double) 0.128e3 - 0.27e2 / (double) 0.128e3 * 
   sqrt(0.9e1 + 0.32e2 * u)));
   nv = ceil(0.4e1 * (double) n * (u - 0.1e1 - vmin)) + 1;
  for (iv = 0; iv <= (int) nv; iv++)
  {v = (double)(vmin + (double) iv * (u - 0.1e1 - vmin) / (double) nv);
  for (ix = 0; ix <= 13 * n - 1; ix++)
  {x = (double)(0.3141592654e1 * (double) ix /(double) n/(double)0.13e2);
  dfv = fabs((0.1e1 + (-0.2e1) * c * c * u) * (cos(0.2e1 * x) + 0.2e1 
  * c * c * v) / (double) (0.1e1 + 0.4e1 * c * c * v * cos(0.2e1 * x) 
  + 0.4e1 *  c * c * c * c * v * v) / (double) sqrt(u + v * cos(0.2e1 
  * x) + (-0.1e1) * c * c * u * u + c * c * v * v) - 0.2e1 * c * 
  sin(0.2e1 * x) / (double) (0.1e1 + 0.4e1 * c * c * v * 
  cos(0.2e1 * x) + 0.4e1 * c * c * c * c * v * v));
  dfx = fabs(0.2e1 * v * (0.1e1 + (-0.2e1) * c * c * u) * 
  sin(0.2e1 * x) / (0.1e1 + 0.4e1 * c * c * v * cos(0.2e1 * x) 
  + 0.4e1 * c * c * c * c * v * v) / (double) sqrt(u + v * cos(0.2e1 
  * x) + (-0.1e1) * c * c * u * u + c * c * v * v) +  0.4e1 * c * v * 
  (cos(0.2e1 * x) + 0.2e1 * c * c * v) / (double) (0.1e1 + 0.4e1 * c 
  * c * v * cos(0.2e1 * x) + 0.4e1 * c * c * c  * c * v * v));
     if ((double) mv < dfv) mv = (double) dfv;
     if ((double) mx < dfx) mx = (double) dfx; } } }
}
printf("mv = %.10f\nmx = %.10f\n", mv , mx );
return(0);}
\end{verbatim}

\bibliography{biblio} \bibliographystyle{alpha}

\begin{thebibliography}{KTJ94}

\bibitem[CL10]{cl}
B.~L. Chalmers and G.~Lewicki.
\newblock A proof of the {G}r\"unbaum conjecture.
\newblock {\em Studia Math.}, 200(2):103--129, 2010.

\bibitem[Gr{\"u}60]{g}
B.~Gr{\"u}nbaum.
\newblock Projection constants.
\newblock {\em Trans. Amer. Math. Soc.}, 95:451--465, 1960.

\bibitem[KS71]{ks}
M.~I. Kadec and M.~G. Snobar.
\newblock Certain functionals on the {M}inkowski compactum.
\newblock {\em Mat. Zametki}, 10:453--457, 1971.

\bibitem[KTJ94]{ktj}
H.~K{\"o}nig and N.~Tomczak-Jaegermann.
\newblock Norms of minimal projections.
\newblock {\em J. Funct. Anal.}, 119(2):253--280, 1994.

\bibitem[Lew88]{l}
D.~R. Lewis.
\newblock An upper bound for the projection constant.
\newblock {\em Proc. Amer. Math. Soc.}, 103(4):1157--1160, 1988.

\end{thebibliography}

\end{document}